\documentclass[12pt]{article}

\usepackage{fontspec, xunicode, xltxtra}  
\usepackage[slantfont,boldfont]{xeCJK}
\usepackage{mathtools}
\usepackage{amsmath,diagbox}
\usepackage{indentfirst}
\usepackage{mathrsfs}
\usepackage{amsfonts}
\usepackage{arydshln}
\usepackage{enumerate}
\usepackage{setspace}
\usepackage{amssymb,amsthm,cases}
\usepackage{amsbsy,pifont}
\usepackage{latexsym,diagbox,tabularx}
\usepackage{amsmath,amscd,bbm,multirow}
\usepackage{graphicx,epsfig,extarrows,dsfont,mathtools}
\usepackage[final,allcolors=blue,colorlinks=true]{hyperref}
\usepackage[normalem]{ulem}
\input xy
\usepackage{caption,subcaption} 
\usepackage{tikz}
\usepackage{tikz-cd}
\usetikzlibrary{decorations.pathreplacing,decorations.markings}
\usetikzlibrary{cd}
\usetikzlibrary{shapes.geometric}
\usetikzlibrary{arrows,arrows.meta}
\usetikzlibrary{graphs,positioning}
\usetikzlibrary{matrix,decorations.pathmorphing}
\usetikzlibrary{math,calc,intersections,through,angles,arrows.meta,shapes.geometric,shadows,quotes,spy,datavisualization,datavisualization.formats.functions,plotmarks}
\tikzset{every picture/.style={samples=300,smooth,line join=round,thick,>=stealth}}
\usepackage{showkeys}
\usepackage{enumerate}
\usepackage{fancyhdr}
\usepackage[sectionbib]{chapterbib}

\tikzset{
	on each segment/.style={
		decorate,
		decoration={
			show path construction,
			moveto code={},
			lineto code={
				\path [#1]
				(\tikzinputsegmentfirst) -- (\tikzinputsegmentlast);
			},
			curveto code={
				\path [#1] (\tikzinputsegmentfirst)
				.. controls
				(\tikzinputsegmentsupporta) and (\tikzinputsegmentsupportb)
				..
				(\tikzinputsegmentlast);
			},
			closepath code={
				\path [#1]
				(\tikzinputsegmentfirst) -- (\tikzinputsegmentlast);
			},
		},
	},
	mid arrow/.style={postaction={decorate,decoration={
				markings,
				mark=at position .5 with {\arrow[#1]{stealth}}
	}}},
}

\numberwithin{equation}{section}

\theoremstyle{plain}
\newtheorem{theorem}{Theorem}[section]

\newtheorem{corollary}[theorem]{Corollary}

\newtheorem{lemma}[theorem]{Lemma}
\newtheorem{proposition}[theorem]{Proposition}

\theoremstyle{definition}
\newtheorem{definition}[theorem]{Definition}

\newtheorem{remark}[theorem]{Remark}

\def\XXint#1#2#3{{\setbox0=\hbox{$#1{#2#3}{\int}$}
		\vcenter{\hbox{$#2#3$}}\kern-.5\wd0}}

\DeclareMathSymbol{\subseteq}{\mathrel}{symbols}{"12}
\DeclareMathSymbol{\supseteq}{\mathrel}{symbols}{"13} 
\DeclareMathSymbol{\subsetneq}{\mathrel}{AMSb}{"28}                  \DeclareMathSymbol{\supsetneq}{\mathrel}{AMSb}{"29}    

\DeclareMathSymbol{\nsubseteq}{\mathrel}{AMSb}{"2A}                  

\DeclareMathSymbol{\nsupseteq}{\mathrel}{AMSb}{"2B}

\DeclareMathDelimiter{\langle}{\mathop}{symbols}{"68}{largesymbols}{"0A}
\DeclareMathDelimiter{\rangle}{\mathclose}{symbols}{"69}{largesymbols}{"0B}

\DeclareSymbolFont{txfontsA}{U}{txmia}{m}{it}
\SetSymbolFont{txfontsA}{bold}{U}{txmia}{bx}{it}
\DeclareFontSubstitution{U}{txmia}{m}{it}
\DeclareMathSymbol{\upalpha}{\mathord}{txfontsA}{"0B}
\DeclareMathSymbol{\upbeta}{\mathord}{txfontsA}{"0C}
\DeclareMathSymbol{\upgamma}{\mathord}{txfontsA}{"0D}
\DeclareMathSymbol{\updelta}{\mathord}{txfontsA}{"0E}
\DeclareMathSymbol{\upepsilon}{\mathord}{txfontsA}{"0F}
\DeclareMathSymbol{\upzeta}{\mathord}{txfontsA}{"10}
\DeclareMathSymbol{\upeta}{\mathord}{txfontsA}{"11}
\DeclareMathSymbol{\uptheta}{\mathord}{txfontsA}{"12}
\DeclareMathSymbol{\upiota}{\mathord}{txfontsA}{"13}
\DeclareMathSymbol{\upkappa}{\mathord}{txfontsA}{"14}
\DeclareMathSymbol{\uplambda}{\mathord}{txfontsA}{"15}
\DeclareMathSymbol{\upmu}{\mathord}{txfontsA}{"16}
\DeclareMathSymbol{\upnu}{\mathord}{txfontsA}{"17}
\DeclareMathSymbol{\upxi}{\mathord}{txfontsA}{"18}
\DeclareMathSymbol{\uppi}{\mathord}{txfontsA}{"19}
\DeclareMathSymbol{\uprho}{\mathord}{txfontsA}{"1A}
\DeclareMathSymbol{\upsigma}{\mathord}{txfontsA}{"1B}
\DeclareMathSymbol{\uptau}{\mathord}{txfontsA}{"1C}
\DeclareMathSymbol{\upupsilon}{\mathord}{txfontsA}{"1D}
\DeclareMathSymbol{\upphi}{\mathord}{txfontsA}{"1E}
\DeclareMathSymbol{\upchi}{\mathord}{txfontsA}{"1F}
\DeclareMathSymbol{\uppsi}{\mathord}{txfontsA}{"20}
\DeclareMathSymbol{\upomega}{\mathord}{txfontsA}{"21}
\DeclareMathSymbol{\upvarepsilon}{\mathord}{txfontsA}{"22}
\DeclareMathSymbol{\upvartheta}{\mathord}{txfontsA}{"23}
\DeclareMathSymbol{\upvarpi}{\mathord}{txfontsA}{"24}
\DeclareMathSymbol{\upvarrho}{\mathord}{txfontsA}{"25}
\DeclareMathSymbol{\upvarsigma}{\mathord}{txfontsA}{"26}
\DeclareMathSymbol{\upvarphi}{\mathord}{txfontsA}{"27}                   

\DeclareSymbolFont{ugmL}{OMX}{mdugm}{m}{n}
\SetSymbolFont{ugmL}{bold}{OMX}{mdugm}{b}{n}
\DeclareMathAccent{\wideparen}{\mathord}{ugmL}{"F3}

\def\pt{\partial}

\def\ra{\rightarrow}

\def\ds{\displaystyle}

\def\s{\subseteq}

\def\e{\epsilon}

\def\ol{\overline}

\def\vp{\varphi}
\def\lg{\langle}
\def\llg{\left\langle}

\def\rg{\rangle}
\def\rrg{\right\rangle}

\def\pt{\partial}

\def\Om{\Omega}
\def\la{\lambda}
\def\al{\alpha}
\def\be{\beta}
\def\de{\delta}
\def\ga{\gamma}

\def\Ga{\Gamma}

\def\ts{\times}

\def\iy{\infty}

\def\f{\frac}

\def\vs{\vspace}

\def\Span{{\rm{span}}}

\def\df{\mathrm d}

\def\wt{\widetilde}
\def\wh{\widehat}

\def\esssup{\operatorname*{ess\ \! sup}}

\def\hra{\hookrightarrow}

\def\diag{{\rm diag}}

\def\mcH{\mathcal{H}}

\def\mcA{\mathcal{A}}

\def\mcV{\mathcal{V}}

\def\mcT{\mathcal{T}}

	\DeclareMathOperator{\Div}{div}

	\DeclareMathOperator{\dist}{dist}

	\DeclareMathOperator{\supp}{{supp}}

	\newcommand{\R}{\mathbb R}
	\newcommand{\N}{\mathbb N}

\begin{document} 
	
	
	\title{A Class of Degenerate Parabolic Equations with Mixed Dirichlet-Neumann Boundary Conditions and Applications}
	
	\author{ Dong-Hui Yang$^{a}$, Mengze Gu$^{b}$ \footnote{\small
		The corresponding author: gu3059268737@gmail.com}, Ghadir Shokor$^{c}$, Yuanzhi Zhou$^{a}$\\
	$^a${\it School of Mathematics and Statistics, Central South University}\\
	{\it Changsha 410075, P.R.China}\\
	$^b${\it School of Mathematics and Science, Changsha Normal University}\\ {\it  Changsha 410083, China}\\
	$^c${\it School of Mathematics and Statistics, Beijing Institute of Technology}\\
	{\it Beijing 100081, P.R.China} }
	
	\date{}
	\maketitle{}
	\thispagestyle{empty}
	\thispagestyle{empty}
	
	\begin{abstract}
		In this paper, we investigate the null controllability of a class of two-dimensional degenerate parabolic equations with mixed Dirichlet--Neumann boundary conditions, where the diffusion coefficient degenerates on a flat portion of the boundary. To overcome the lack of higher regularity near the degenerate boundary when performing integration by parts, we introduce a shape design method tailored to mixed boundary value problems. Combined with a family of uniformly parabolic approximate equations, this approach allows us to establish uniform global Carleman estimates and subsequently derive the observability inequality for the original degenerate system. 
		
		\medskip
		
		
		
	\end{abstract}

\section{Introduction}\label{S1}

The controllability of uniformly  parabolic partial differential equations (PDEs) has been extensively developed since the 1960s and has been investigated in numerous works \cite{Coron,Fursikov,Komornik,Lions,Lu,Miller,Privat,Rousseau,Weiss,Yang5,Zuazua}. In contrast, the controllability of degenerate PDEs has been studied since the 2000s, with a growing body of literature \cite{Alabau,Araruna,Barcena,Beauchard,Buffe,Cannarsa1,Cannarsa,Guo,Yang,Yang7,Yang8,Yang4,Yang9,Yang5,Yang1,Yang2,Yang3,Yang6,Zhang}. Prior to 2020, most of these results were restricted to the one-dimensional case, with only a few addressing higher-dimensional settings \cite{Araruna,Beauchard,Cannarsa}. According to J.-L. Lions' Hilbert Uniqueness Method (HUM), controllability is equivalent to observability, a principle that also holds for degenerate PDEs (see \cite{Lions,Yang1}).

Shape design (or shape optimization), a well-established branch of mathematics with broad applications in industry and engineering, has been studied in many works \cite{Chenais,Greco,Guo2,Guo1,Henrot,Privat,Wang,Wang1}. Recently, the authors applied shape design techniques to investigate the controllability of degenerate PDEs \cite{Guo,Yang8,Yang3,Yang6}. It has been demonstrated that shape design is a powerful approach for studying the controllability of degenerate PDEs, as it avoids the need for higher regularity on the degenerate part—a crucial advantage when performing integration by parts.

In this paper, we consider the following degenerate parabolic equation:
\begin{equation}\label{11.22.1}
	\begin{cases}
		\partial_{t}\vp -\Div(A \nabla \vp)=f, &\text{in }Q, \\ 
		\vp=0, &\text{on } \Ga\ts (0,T), \\ 
		\f{\pt \vp}{\pt \nu_A}=0, &\text{on }\Ga_2^0\ts (0,T),\\
		\vp(0)=\vp^0, &\text{in }\Om,
	\end{cases}
\end{equation}
where $A=\diag(1,x_2^\al)$ with $\al\in [1,2)$ being a given constant, $\Om\s \R_+^2=\{x=(x_1,x_2)\colon x_2>0\}$ is a bounded convex domain with boundary $\pt\Om\in C^3$ such that
\begin{equation*}
	\Ga_2^0=\{x\in\ol\Om\colon x_2=0\}=[a_0,b_0]\ts\{0\}  \text{ and }   \Ga=\pt\Om\setminus\Ga_2^0,
\end{equation*}
where $a_0<0<b_0$ (so that $(0,0)\in\R^2$ is a relative interior point of $\Ga_2^0$ with respect to $\pt\Om$), $Q=\Om\ts (0,T)$ for a given constant $T>0$, $\vp^0\in L^2(\Om)$ is the initial datum, and $f\in L^2(Q)$ is the nonhomogeneous term. The conormal derivative $\f{\pt\vp}{\pt \nu_A}=A\nabla\vp\cdot\nu$ (where $\nu$ denotes the outward unit normal vector on $\pt\Om$) will be rigorously defined in Section \ref{S2}. 

We denote 
\begin{equation}\label{11.22.2}
	\mcA \vp=-\Div(A\nabla \vp),\quad M=\sup_{x\in \Om} |x|, 
\end{equation}
and fix $\de_0>0$ sufficiently small such that there exists a convex function $h\in C^3([a,b])$ (with $a<a_0<0<b_0<b$, $h(a)=h(b)=8\de_0$, and $h\equiv 0$ on $[a_0,b_0]$) satisfying 
\begin{equation}\label{11.22.3}
	\pt\Om\cap \{0<x_2<8\de_0\}=\left\{(x_1,x_2)\in\R^2\colon x_1\in (a,b) \text{ and } x_2=h(x_1)\in(0,8\de_0)\right\}
\end{equation}
(equivalently, $\Om\cap\{0<x_2<8\de_0\}=\{(x_1,x_2)\in\R^2\colon x_1\in(a,b) \text{ and } h(x_1)<x_2<8\de_0\}$).

In \cite{Yang8}, the shape design approach yields boundary observability in the weakly degenerate Dirichlet setting, while its extension to Neumann boundary conditions is left open. The present work addresses this issue in the strongly degenerate regime $1\leq\alpha<2$, with a homogeneous conormal Neumann condition on the degenerate boundary portion $\Gamma_2^0$ and a Dirichlet condition on the remaining boundary. To the best of our knowledge, this is the first work to establish boundary observability, and hence the corresponding null controllability by duality, for problem~\eqref{11.22.1} on the class of smooth convex planar domains considered here, without assuming a Cartesian product structure. This is achieved through a shape design adapted to the mixed boundary conditions and Carleman estimates that are uniform with respect to the approximation parameter.

The rest of this paper is organized as follows. In Section \ref{S2}, we present some preliminary results, including weighted Sobolev spaces, a weighted divergence space, and the corresponding Green's formulas. In Section \ref{S3}, we introduce a shape design method for mixed Dirichlet-Neumann boundary value problems and establish a crucial result, Theorem \ref{11.23.T3}. In Section \ref{S4}, we employ a family of uniformly parabolic approximate equations to derive uniform Carleman estimates, which in turn yield the observability inequality for the original degenerate equation. Finally, in Section \ref{S5}, we provide the proofs of two fundamental results.

The main results of this paper are Theorems \ref{11.23.T1}, \ref{11.23.T3}, and \ref{11.26.T1}.

\section{Preliminaries}\label{S2}

In this section, we present some preliminary results, including weighted Sobolev spaces, a weighted divergence space, and the corresponding Green's formulas.

\subsection{Weighted Sobolev spaces}\label{S2.1}

Define
\begin{equation*}
	\wh H^1(\Om;A)=\left\{u\in L^2(\Om)\colon \int_\Om \nabla u\cdot A\nabla u\df x<+\iy\right\}, 
\end{equation*}
where $\nabla u=(\pt_{x_1}u,\pt_{x_2}u)$ denotes the distributional gradient of $u$. Its inner product and norm are defined by 
\begin{equation*}
	(u,v)_{\wh H^1(\Om;A)}=\int_\Om \left(uv+\nabla u\cdot A\nabla v\right)\df x,\quad \|u\|_{\wh H^1(\Om;A)}=(u,u)_{\wh H^1(\Om;A)}^\f{1}{2}. 
\end{equation*}
Set
\begin{equation*}
	H^1(\Om;A)=\text{the closure of } C^\iy(\ol\Om) \text{ in } \wh H^1(\Om;A). 
\end{equation*}

Define $\dist(X,Y)=\inf_{x\in X, y\in Y}|x-y|$ for non-empty subsets $X,Y\s\R^2$, and 
\begin{equation*}
	C_\Ga^\iy(\ol\Om)=\left\{u\in C^\iy(\ol\Om)\colon \dist(\supp u, \ol \Ga)>0\right\}. 
\end{equation*}
Set
\begin{equation*}
	\begin{split} 
		H_\Ga^1(\Om;A)
		&=\text{the closure of } C_\Ga^\iy(\ol\Om) \text{ in } H^1(\Om;A), \\
		H_\Ga^{-1}(\Om;A)
		&=\text{the dual space of } H_\Ga^1(\Om;A) \text{ with respect to the pivot space } L^2(\Om). 
	\end{split} 
\end{equation*}

Define
\begin{equation*}
	H^2(\Om;A)=\left\{u\in H^1(\Om;A)\colon \mcA u\in L^2(\Om)\right\}.
\end{equation*}
The inner product and norm in $H^2(\Om;A)$ are defined by 
\begin{equation*}
	(u,v)_{H^2(\Om;A)}=\int_\Om \left(uv+\nabla u\cdot A\nabla v+(\mcA u)(\mcA v)\right)\df x, \quad \|u\|_{H^2(\Om;A)}=(u,u)_{H^2(\Om;A)}^\f{1}{2}. 
\end{equation*}
Set
\begin{equation*}
	D(\mcA)=\{u\in H_\Ga^1(\Om;A)\colon \mcA u\in L^2(\Om)\}. 
\end{equation*}

It is easily verified that 
\begin{equation*}
	H_\Ga^1(\Om;A),\quad H_\Ga^{-1}(\Om;A), \quad  H^1(\Om;A),\quad \wh H^1(\Om;A), \quad H^2(\Om;A),\quad  D(\mcA)
\end{equation*}
are Hilbert spaces (see \cite{GC,Heinonen}), where the inner product on $D(\mcA)$ is defined by 
\begin{equation*}
	(u,v)_{D(\mcA)}=\int_\Om \left( uv+\nabla u\cdot A\nabla v+(\mcA u)(\mcA v)\right)\df x. 
\end{equation*}

\begin{lemma}\label{11.22.L9}
	It holds that $H^1(\Om;A)=\wh H^1(\Om;A)$. 
\end{lemma}

\begin{proof}
	By definition, $H^1(\Om;A)\s \wh H^1(\Om;A)$, so it suffices to prove that $C^\iy(\ol\Om)$ is dense in $\wh H^1(\Om;A)$. Let $v\in \wh H^1(\Om;A)$ be arbitrary.

	Recall the boundary representation \eqref{11.22.3}, with
	$h(a)=h(b)=8\delta_0$. If $0<\delta<\delta_0$ and
	$(x_1,x_2)\in\Omega\cap\{0<x_2<\delta\}$, then
	\[
	h(x_1)<x_2<\delta<2\delta-x_2<2\delta<2\delta_0<8\delta_0.
	\]
	Thus $(x_1,2\delta-x_2)\in\Omega\cap\{\delta<x_2<2\delta\}$,
	and the even reflection used below is well defined.
	\begin{equation*}
		v_\de(x_1,x_2) =
		\begin{cases}
			v(x_1,x_2), & \text{if } (x_1,x_2)\in \wh\Om_\de := \Om\cap\{x_2\geq \de\},\\
			v(x_1, 2\de-x_2), & \text{if } (x_1,x_2)\in \Om-\wh\Om_\de = \Om\cap\{0<x_2<\de\}.
		\end{cases}
	\end{equation*}
	
	Since $x_2^\al\geq \de^\al>0$ on $\wh\Om_\de$, we have $v|_{\wh\Om_\de}\in H^1(\wh\Om_\de)$. By standard properties of even reflection across $\{x_2=\de\}$, it follows that $v_\de\in H^1(\Om)$. Moreover, on $\Om-\wh\Om_\de$, setting $y_2 = 2\de-x_2 > x_2 > 0$, we have $x_2^\al < y_2^\al$ and $|\nabla v_\de(x_1,x_2)| = |(\pt_{x_1}v(x_1,y_2), -\pt_{x_2}v(x_1,y_2))|$, which yields
	\begin{equation*}
		\begin{split} 
			\|v_\de - v\|_{\wh H^1(\Om;A)}^2
			&= \int_{\Om-\wh\Om_\de} \left(|v_\de - v|^2 + \nabla(v_\de - v)\cdot A\nabla(v_\de - v)\right)\df x\\
			&\leq 2\int_{\Om-\wh\Om_\de} \left(|v|^2 + \nabla v\cdot A\nabla v\right)\df x + 2\int_{\Om\cap\{\de<y_2<2\de\}} \left(|v|^2 + \nabla v\cdot A\nabla v\right)\df y.
		\end{split}
	\end{equation*}
	By the absolute continuity of the Lebesgue integral, the right-hand side tends to $0$ as $\de\to 0^+$. Finally, for each fixed $\de\in (0,\de_0 )$, since $v_\de\in H^1(\Om)$ and $C^\iy(\ol\Om)$ is dense in $H^1(\Om)$ (with $H^1(\Om)\hra \wh H^1(\Om;A)$ continuously), we can choose $\phi_\de\in C^\iy(\ol\Om)$ such that $\|\phi_\de - v_\de\|_{\wh H^1(\Om;A)} < \de$. Therefore, $\phi_\de \to v$ in $\wh H^1(\Om;A)$ as $\de\to 0^+$, completing the proof.
\end{proof}

\begin{lemma}[Hardy-like inequalities]\label{01.28.L3}
	Let $u \in H_\Ga^1(\Om;A)$.
	
	{\rm (i)} If $\alpha \in (1,2)$, then there holds
	\begin{equation}\label{02.08.1}
		\int_{\Omega} x_2^{\alpha - 2}u^2 \, \df x \leq \frac{4}{(\alpha - 1)^2}\int_{\Omega} x_2^{\alpha} \left|\frac{\partial u}{\partial x_2}\right|^2 \, \df x.
	\end{equation}
	
	{\rm (ii)} If $\alpha = 1$, then there holds
	\begin{equation}\label{02.08.1b}
		\int_\Omega u^2\, \df x\leq 4M\int_\Omega x_2\left|\frac{\partial u}{\partial x_2}\right|^2\, \df x.
	\end{equation}
	Moreover, for any $\beta>0$ and $\lambda>0$, there exists a constant $C(\beta,\lambda)>0$ depending only on $\beta$, $\lambda$, and $M$ such that
	\begin{equation}\label{02.08.2}
		\begin{split} 
			\lambda\int_\Omega x_2^{-1+\beta}u^2\, \df x
			&\leq \int_\Omega x_2 \left|\frac{\partial u}{\partial x_2}\right|^2\, \df x+C(\beta,\lambda)\int_\Omega u^2\, \df x\\
			&\leq \big(1+4MC(\beta,\lambda)\big)\int_\Omega x_2\left|\frac{\partial u}{\partial x_2}\right|^2\, \df x.
		\end{split} 
	\end{equation} 
\end{lemma}

\begin{proof}
	{\rm (i)} Let $\alpha \in (1,2)$. By the definition of $H_\Ga^1(\Om;A)$ and Fatou's lemma, it suffices to establish \eqref{02.08.1} for $u \in C_\Ga^\iy(\ol\Om)$. Extending $u$ by zero outside $\Om$, we have $u \in C^\iy(\R\ts [0,+\iy))$ with $\supp u \s [-M,M]\ts [0,M)$, where $M$ is defined in \eqref{11.22.2}. In particular, $u(x_1,M)=0$ for all $x_1\in (-M,M)$. Applying H\"older's inequality gives
	\begin{equation*}
		\begin{split}  
			|u(x_1,x_2)|^2
			&= \left|\int_{x_2}^M\frac{\partial u}{\partial x_2} (x_1,s) \, \df s\right|^2 \leq \left(\int_{x_2}^Ms^{\frac{\alpha + 1}{2}}\left|\frac{\partial u}{\partial x_2}(x_1,s)\right|^2 \, \df s\right)\left(\int_{x_2}^M s^{-\frac{\alpha + 1}{2}} \, \df s\right)\\
			&\leq \frac{2}{\alpha - 1} x_2^{-\frac{\alpha - 1}{2}} \int_{x_2}^Ms^{\frac{\alpha + 1}{2}}\left|\frac{\partial u}{\partial x_2}(x_1,s)\right|^2 \, \df s.
		\end{split}
	\end{equation*}
	Multiplying both sides by $x_2^{\alpha - 2}$, integrating with respect to $x_2$ over $(0,M)$, and changing the order of integration via Fubini's theorem, we obtain
	\begin{equation*}
		\begin{split}
			\int_0^M x_2^{\alpha - 2}|u(x_1,x_2)|^2 \, \df x_2
			&\leq \frac{2}{\alpha - 1}\int_0^M\left(\int_{x_2}^M s^{\frac{\alpha + 1}{2}}\left|\frac{\partial u}{\partial x_2}(x_1,s)\right|^2 \, \df s\right)x_2^{\frac{\alpha - 3}{2}} \, \df x_2\\
			&=\frac{4}{(\alpha - 1)^2}\int_0^M s^{\alpha} \left|\frac{\partial u}{\partial x_2}(x_1,s)\right|^2 \, \df s.
		\end{split}
	\end{equation*}
	Integrating with respect to $x_1 \in (-M,M)$ yields \eqref{02.08.1}.
	
	{\rm (ii)} Now let $\alpha = 1$ and $u\in H_\Ga^1(\Om;A)$. Since $0 < x_2 < M$ on $\Omega$, we have $\int_\Omega x_2 u^2\, \df x\leq M\int_\Omega u^2\, \df x$ and
	\begin{equation*}
		\begin{split}
			\int_\Omega |\nabla (x_2 u^2)|\, \df x
			&=2\int_\Omega |u|\left|\left(x_2\frac{\partial u}{\partial x_1},\, \f{1}{2}u+x_2\frac{\partial u}{\partial x_2}\right)\right|\, \df x\\
			&\leq C\int_\Omega u^2\, \df x+C\int_\Omega A\nabla u\cdot \nabla u\, \df x<+\infty, 
		\end{split}
	\end{equation*}
	where $C>0$ depends only on $M$. 
	Furthermore, by the density of $C_\Ga^\iy(\ol\Om)$ in $H_\Ga^1(\Om;A)$ and the vanishing weight $(x_2u^2)|_{x_2=0}=0$, we have $x_2 u^2\in W_0^{1,1}(\Omega)$, and hence $\int_\Omega \frac{\partial}{\partial x_2}(x_2 u^2)\,\df x = 0$. Using Young's inequality, we deduce that
	\begin{equation*}
		\begin{split}
			\int_\Omega u^2\, \df x
			&=\int_\Omega \frac{\partial x_2}{\partial x_2}u^2\, \df x=\int_\Omega \frac{\partial}{\partial x_2}\left(x_2 u^2\right)\, \df x-2\int_\Omega x_2 u\frac{\partial u}{\partial x_2}\, \df x\\
			&\leq \frac{1}{2M}\int_\Omega x_2 u^2\, \df x + 2M\int_\Omega x_2\left|\frac{\partial u}{\partial x_2}\right|^2\, \df x \leq \frac{1}{2}\int_\Omega u^2\, \df x + 2M\int_\Omega x_2\left|\frac{\partial u}{\partial x_2}\right|^2\, \df x,
		\end{split}
	\end{equation*}
	which implies \eqref{02.08.1b}.
	
	Finally, for any $\beta>0$ and $\lambda>0$, set
	\begin{equation*}
		\begin{split} 
		d 
		&= \min\left\{\frac{M}{2}, \left(\frac{\beta}{2(\lambda + 1)}\right)^{\frac{2}{\beta}}\right\}, \\ 
		C(\beta, \lambda) 
		&= \max\left\{1,\sup_{x_2 \in (d, M)} \left[\left((\lambda + 1)x_2^{\frac{\beta}{2}} - \frac{\beta}{2}\right)x_2^{-1 + \frac{\beta}{2}}\right] \right\}< +\infty.
		\end{split} 
	\end{equation*}
	For each $\delta\in (0,d)$, denote $\wh \Omega_\delta = \Om\cap \{x_2>\de\}$ and $\Gamma_2^\delta = \Om\cap\{x_2=\delta\}$. Since $u \in H^1(\wh\Omega_\delta)$ with $u=0$ on $\partial\wh\Omega_\delta - \Gamma_2^\delta$, expanding the nonnegative integral
	\begin{equation*}
		0\leq\int_{\wh\Omega_\delta}\left(x_2^{\frac{1}{2}}\frac{\partial u}{\partial x_2}+x_2^{\frac{-1+\beta}{2}}u\right)^2\, \df x
		=\int_{\wh\Omega_\delta} x_2 \left(\frac{\partial u}{\partial x_2}\right)^2\, \df x+\int_{\wh\Omega_\delta} x_2^{-1+\beta} u^2\, \df x+2\int_{\wh\Omega_\delta} x_2^{\frac{\beta}{2}}u\frac{\partial u}{\partial x_2}\, \df x
	\end{equation*}
	and integrating the last term by parts with respect to $x_2$ gives
	\begin{equation*}
		\begin{split}
			0
			&\leq \int_{\wh\Omega_\delta}x_2\left(\frac{\partial u}{\partial x_2}\right)^2\, \df x+\int_{\wh\Omega_\delta}x_2^{-1+\beta}u^2\, \df x -\int_{\Gamma_2^\delta}x_2^{\frac{\beta}{2}}u^2\, \df x_1-\frac{\beta}{2}\int_{\wh\Omega_\delta} x_2^{-1+\frac{\beta}{2}}u^2\, \df x.
		\end{split}
	\end{equation*}
	Rearranging terms, we obtain
	\begin{equation*}
		\begin{split}
			&\int_{\Gamma_2^\delta}x_2^{\frac{\beta}{2}}u^2\, \df x_1+\lambda\int_{\wh\Omega_\delta} x_2^{-1+\beta}u^2\, \df x\\
			&\leq \int_{\wh\Omega_\delta} x_2\left(\frac{\partial u}{\partial x_2}\right)^2\, \df x+(\lambda+1)\int_{\wh\Omega_\delta} x_2^{-1+\beta}u^2\, \df x-\frac{\beta}{2}\int_{\wh\Omega_\delta} x_2^{-1+\frac{\beta}{2}}u^2\, \df x\\
			&= \int_{\wh\Omega_\delta} x_2\left(\frac{\partial u}{\partial x_2}\right)^2\, \df x+\int_{\wh\Omega_\delta} \left((\lambda+1)x_2^{\frac{\beta}{2}}-\frac{\beta}{2}\right)x_2^{-1+\frac{\beta}{2}}u^2\, \df x.
		\end{split}
	\end{equation*}
	Note that $(\lambda+1)x_2^{\frac{\beta}{2}}-\frac{\beta}{2}\leq 0$ for all $x_2\in [0, (\frac{\beta}{2(\lambda+1)})^{\frac{2}{\beta}}]$. Dropping the nonnegative boundary integral on $\Gamma_2^\delta$ and splitting the domain of integration into $\wh\Omega_\delta \cap \{x_2 \leq d\}$ and $\wh\Omega_\delta \cap \{x_2 > d\}$, we get
	\begin{equation*}
		\lambda \int_{\wh\Omega_\delta} x_2^{-1 + \beta} u^2 \, \df x \leq \int_{\wh\Omega_\delta} x_2 \left(\frac{\partial u}{\partial x_2}\right)^2 \, \df x + C(\beta, \lambda) \int_{\Om} u^2 \, \df x.
	\end{equation*}
	Passing to the limit as $\delta \to 0^+$ via the monotone convergence theorem and combining with \eqref{02.08.1b}, we arrive at \eqref{02.08.2}. This completes the proof.
\end{proof}

\begin{remark}\label{11.22.R1}
	By Lemma \ref{01.28.L3} (together with $0 < x_2 < M$ on $\Om$), for any $u\in H_\Ga^1(\Om;A)$, there exists a positive constant $C_P$, depending only on $\al$ and $M$, such that 
	\begin{equation}\label{11.22.4}
		\int_\Om u^2\,\df x\leq C_P\int_\Om \nabla u\cdot A\nabla u\,\df x,
	\end{equation}
	which is the Poincar\'e inequality in $H_\Ga^1(\Om;A)$. Throughout this paper, we equip $H_\Ga^1(\Om;A)$ with the equivalent norm
	\begin{equation*}
		\|u\|_{H_\Ga^1(\Om;A)}=\left(\int_\Om \nabla u\cdot A\nabla u\,\df x\right)^\f{1}{2}.
	\end{equation*}
\end{remark}

\begin{lemma}\label{11.22.L1}
	The embedding $H_\Ga^1(\Om;A)\hra L^2(\Om)$ is compact. 
\end{lemma}

\begin{proof}
	Let $\{u_n\}_{n=1}^\iy \s H_\Ga^1(\Om;A)$ be a bounded sequence with $\|u_n\|_{H_\Ga^1(\Om;A)}\leq K$ for all $n\geq 1$. By the reflexivity of $H_\Ga^1(\Om;A)$, up to a subsequence, we may assume that $u_n \ra u$ weakly in $H_\Ga^1(\Om;A)$ for some $u\in H_\Ga^1(\Om;A)$. Set $v_n = u_n - u$, so that $v_n \ra 0$ weakly in $H_\Ga^1(\Om;A)$ and $\|v_n\|_{H_\Ga^1(\Om;A)}\leq 2K$.
	
	For any $\de>0$, denote $\wh\Om_\de=\Om\cap \{x_2>\de\}$. By Lemma \ref{01.28.L3}, we have
	\begin{equation*}
		\int_{\Om- \wh\Om_\de} v_n^2\,\df x \leq
		\begin{cases}
			\de^{2-\al}\ds\int_\Om x_2^{\al-2}v_n^2\,\df x \leq C\de^{2-\al}\|v_n\|_{H_\Ga^1(\Om;A)}^2 \leq C\de^{2-\al}K^2, & \text{if }\al\in (1,2),\\[3ex]
			\de^{\frac{1}{2}}\displaystyle\int_\Om x_2^{-\frac{1}{2}}v_n^2\,\df x \leq C\de^{\frac{1}{2}}\|v_n\|_{H_\Ga^1(\Om;A)}^2 \leq C\de^{\frac{1}{2}}K^2, & \text{if }\al=1,
		\end{cases}
	\end{equation*}
	where $C>0$ is independent of $n$ and $\de$. Thus, for any $\e>0$, we can fix $\de>0$ sufficiently small such that 
	\begin{equation*}
		\int_{\Om- \wh\Om_\de} v_n^2\,\df x < \frac{\e}{2}, \text{ for all } n\geq 1.
	\end{equation*}
	
	On the other hand, since $x_2^\al\geq \de^\al>0$ on $\wh\Om_\de$, the sequence $\{v_n|_{\wh\Om_\de}\}$ is bounded in $H^1(\wh\Om_\de)$ and converges weakly to $0$ in $H^1(\wh\Om_\de)$. By the classical Rellich-Kondrachov theorem on the Lipschitz domain $\wh\Om_\de$, $v_n|_{\wh\Om_\de} \to 0$ strongly in $L^2(\wh\Om_\de)$ as $n\to\iy$. Hence, there exists $N\geq 1$ such that $\int_{\wh\Om_\de} v_n^2\,\df x < \frac{\e}{2}$ for all $n\geq N$, which yields $\|v_n\|_{L^2(\Om)}^2 < \e$ for all $n\geq N$. Therefore, $u_n \to u$ strongly in $L^2(\Om)$.
\end{proof}

\begin{lemma}\label{11.24.L1}
	Let $\al\in [1,2)$. Then $C_0^\iy(\Om)$ is dense in $H_\Ga^1(\Om;A)$. Moreover,  
	\begin{equation}\label{11.24.3}
		H_{\Ga}^1(\Om;A)=\left\{u\in H^1(\Om;A)\colon Tu=0 \mbox{ on }\Ga\right\}.
	\end{equation}
	Here $T$ denotes the local Sobolev trace on the nondegenerate
	boundary $\Gamma$. More precisely, $Tu=0$ on $\Gamma$ means
	that, for every relatively compact boundary patch
	$\Gamma'\Subset\Gamma$, the classical trace of $u$ on
	$\Gamma'$ vanishes. This trace is well defined because
	$u\in H^1$ in a neighborhood of $\Gamma'$, where $A$ is
	uniformly elliptic. No classical trace on $\Gamma_2^0$
	is used in \eqref{11.24.3}.
	We also write
	\[
	H_0^1(\Omega;A)
	:=\overline{C_c^\infty(\Omega)}^{\,H^1(\Omega;A)}
	=H_\Gamma^1(\Omega;A),
	\]
	where the equality follows from the density assertion of
	this lemma. On this common space we use the energy norm
	specified in Remark~\ref{11.22.R1}.
\end{lemma}

\begin{proof}
	We denote the right-hand side of \eqref{11.24.3} by $\mcH$. 
	
For every $\Gamma'\Subset\Gamma$, the local classical trace is continuous with respect to the $H^1(\Omega;A)$ norm. Every function in $C_\Gamma^\infty(\overline\Omega)$ has zero trace on $\Gamma'$. Passing to the closure, and then varying $\Gamma'$, proves $H_\Gamma^1(\Omega;A)\subset\mathcal H$.

Conversely, let $u\in \mcH$. For every $\e>0$, there exists $N\in\N$ such that $\|u_N-u\|_{H^1(\Om;A)}<\f{1}{2}\e$, where 
	\begin{equation*}
		u_N=\max\{-N, \min\{u, N\}\}. 
	\end{equation*}
	Moreover, $Tu_N=0$ on $\Ga$. 
	
	Let $\de\in (0,1)$. When $\al>1$, we take the Lipschitz cut-off function 
	\begin{equation*}
		\zeta(x_2)=
		\begin{cases}
			0,& x_2\leq \de, \\
			\de^{-1}(x_2-\de), &x_2\in (\de,2\de), \\
			1, &x_2\geq 2\de. 
		\end{cases}
	\end{equation*}
	Then 
	\begin{equation*}
		\int_\Om x_2^\al |\zeta'(x_2)|^2\df x\leq \f{2M}{\de^2} \int_\de^{2\de} x_2^\al \df x_2\leq 6M\de^{\al-1}\ra 0\mbox{ as }\de\ra 0^+.  
	\end{equation*}
	When $\al=1$, we take the logarithmic cut-off function 
	\begin{equation*}
		\zeta(x_2)=
		\begin{cases}
			0, &x_2\leq \de^2, \\
			\f{\log\f{x_2}{\de^2}}{\log\f{1}{\de}}, &x_2\in (\de^2,\de),\\
			1, &x_2\geq \de.
		\end{cases}
	\end{equation*}
	Then 
	\begin{equation*}
		\int_\Om x_2|\zeta'(x_2)|^2\df x\leq \f{2M}{\log^2\f{1}{\de}}\int_{\de^2}^\de \f{\df x_2}{x_2} =\f{2M}{|\log\de|}\ra 0\mbox{ as } \de\ra 0^+. 
	\end{equation*}
	
	For fixed $N$, the full weighted $H^1$ norm satisfies
	\[
	\begin{split}
		\|u_N-\zeta u_N\|_{H^1(\Omega;A)}^2
		&\leq \int_\Omega(1-\zeta)^2|u_N|^2\,\mathrm dx +2\int_\Omega(1-\zeta)^2
		\nabla u_N\cdot A\nabla u_N\,\mathrm dx\\
		&\quad+2\int_\Omega x_2^\alpha|\zeta'|^2|u_N|^2\,\mathrm dx
		\longrightarrow0\qquad(\delta\to0^+).
	\end{split}
	\]
	Indeed, the first two integrals tend to zero by dominated
	convergence, and the third is bounded by
	$2N^2\int_\Omega x_2^\alpha|\zeta'|^2\,\mathrm dx\to0$.
	
	Since $\zeta(x_2)=0$ near $\{x_2=0\}$ and $Tu_N=0$ on $\Ga$, we have $\zeta u_N\in H_0^1(\Om)$. Choosing $\de>0$ sufficiently small such that $\|\zeta u_N - u_N\|_{H^1(\Om;A)}<\f{1}{2}\e$, we obtain $\|\zeta u_N-u\|_{H^1(\Om;A)}<\e$. Since $H_0^1(\Om)\hra H_\Ga^1(\Om;A)$ and $C_0^\iy(\Om)$ is dense in $H_0^1(\Om)$, it follows that $C_0^\iy(\Om)$ is dense in $H_0^1(\Om;A)$, and $\mcH\s H_\Ga^1(\Om;A)$. This proves the lemma. 
\end{proof}

Denote 
\begin{equation*}
	H_\Ga^1(\Om)=\text{the closure of } C_\Ga^\iy(\ol\Om) \text{ in } H^1(\Om). 
\end{equation*}
Since $0 < x_2 < M$ on $\Om$, we have $H^1(\Om)\hra H^1(\Om;A)$ continuously, and consequently $H_\Ga^1(\Om)\hra H_\Ga^1(\Om;A)$ continuously. Moreover, $H_\Ga^1(\Om)$ is dense in $H_\Ga^1(\Om;A)$ because $C_\Ga^\iy(\ol\Om)\s H_\Ga^1(\Om)$ is dense in $H_\Ga^1(\Om;A)$.

We introduce the classical trace operator
\begin{equation}\label{11.22.5}
	\ga\colon H_\Ga^1(\Om)\ra H_{00}^\f{1}{2}(\Ga_2^0),
\end{equation}
obtained by extending by continuity the mapping
\begin{equation*}
	\ga u=u|_{\Ga_2^0}, \text{ for all } u\in C_\Ga^\iy(\ol\Om), 
\end{equation*}
where $H_{00}^\f{1}{2}(\Ga_2^0)$ denotes the Lions-Magenes space (see \cite{Lions}) of functions in $H^\f{1}{2}(\Ga_2^0)$ whose extension by zero to $\R\ts \{0\}$ belongs to $H^\f{1}{2}(\R\ts\{0\})$, equipped with the norm
\begin{equation*}
	\|\phi\|_{H_{00}^\f{1}{2}(\Ga_2^0)} = \|\wt\phi\|_{H^\f{1}{2}(\R\ts\{0\})},
\end{equation*}
where $\wt\phi$ is the zero extension of $\phi$ to $\R\ts\{0\}$.

\begin{lemma}\label{11.22.L2}
	The trace operator $\ga$ defined in \eqref{11.22.5} is a bounded linear operator. 
\end{lemma}

\begin{proof}
	By the density of $C_\Ga^\iy(\ol\Om)$ in $H_\Ga^1(\Om)$ and linearity, it suffices to prove that there exists a constant $C>0$, depending only on $\Om$, such that 
	\begin{equation}\label{11.22.L2.est}
		\|\ga u\|_{H_{00}^\f{1}{2}(\Ga_2^0)} \leq C\|u\|_{H^1(\Om)}, \text{ for all } u\in C_\Ga^\iy(\ol\Om).
	\end{equation}
	
	Let $u\in C_\Ga^\iy(\ol\Om)$. Since $\Om\s \R_+^2$ is a bounded domain with $\pt\Om\in C^3$, by the classical Sobolev trace theorem (see, e.g., \cite{Lions}), the standard boundary trace $u|_{\pt\Om}$ belongs to $H^\f{1}{2}(\pt\Om)$ and satisfies
	\begin{equation*}
		\|u|_{\pt\Om}\|_{H^\f{1}{2}(\pt\Om)} \leq C_1\|u\|_{H^1(\Om)},
	\end{equation*}
	where $C_1>0$ depends only on $\Om$. Moreover, because $\dist(\supp u, \ol\Ga)>0$, we have $u|_{\pt\Om}=0$ on $\Ga=\pt\Om-\Ga_2^0$, so the zero extension $\wt{\ga u}$ of $\ga u=u|_{\Ga_2^0}$ to $\R\ts\{0\}$ coincides with the zero extension of $u|_{\pt\Om}$ outside $\pt\Om$ when viewed in local boundary coordinates flattening $\pt\Om$ near $\ol{\Ga_2^0}$. Specifically, since $\pt\Om\in C^3$ and $\Ga_2^0=(a_0,b_0)\ts\{0\}\s \pt\Om$ is flat, we can choose an open neighborhood $U\s \R^2$ of $\ol{\Ga_2^0}$ and a $C^3$-diffeomorphism $\Phi\colon U\ra V\s \R^2$ such that $\Phi(x_1,x_2)=(x_1,x_2-h(x_1))$ for all $(x_1,x_2)\in U\cap(\R\ts\{0\})$ and $\Phi(U\cap\pt\Om)=V\cap(\R\ts\{0\})$. Since $\supp(u|_{\pt\Om})\s \Ga_2^0$ is a compact subset of $U\cap\pt\Om$, transporting $u|_{\pt\Om}$ via $\Phi^{-1}$ preserves the $H^\f{1}{2}$-norm up to a constant depending only on $\Om$, yielding
	\begin{equation*}
		\|\wt{\ga u}\|_{H^\f{1}{2}(\R\ts\{0\})} \leq C_2\|u|_{\pt\Om}\|_{H^\f{1}{2}(\pt\Om)} \leq C_1C_2\|u\|_{H^1(\Om)},
	\end{equation*}
	which gives \eqref{11.22.L2.est}. Hence, $\ga$ extends uniquely to a bounded linear operator from $H_\Ga^1(\Om)$ into $H_{00}^\f{1}{2}(\Ga_2^0)$.
\end{proof}

\begin{remark}\label{11.22.R2}
	By the definition of $H_{00}^\f{1}{2}(\Ga_2^0)$ (see \cite{Lions}), the space $C_0^\iy(\Ga_2^0)$ is dense in $H_{00}^\f{1}{2}(\Ga_2^0)$. 
\end{remark}

\subsection{Weighted divergence spaces}\label{S2.2}

Define
\begin{equation*}
	H_{\Div}(\Om)=\left\{\mathbf{u}\in [L^2(\Om)]^2\colon \Div\mathbf{u}\in L^2(\Om)\right\},
\end{equation*}
where $\Div\mathbf{u}$ denotes the distributional divergence of $\mathbf{u}$. Its inner product and norm are defined by 
\begin{equation*}
	(\mathbf{u},\mathbf{v})_{H_{\Div}(\Om)}=\int_\Om \left(\mathbf{u}\cdot \mathbf{v}+(\Div\mathbf{u})(\Div\mathbf{v})\right)\,\df x,\quad \|\mathbf{u}\|_{H_{\Div}(\Om)}=(\mathbf{u},\mathbf{u})_{H_{\Div}(\Om)}^\f{1}{2}.  
\end{equation*}
We also introduce the weighted divergence space
\begin{equation*}
	H_{\Div}(\Om;A^{-1})=\left\{\mathbf{u}\in [L^2(\Om)]^2\colon \mathbf{u}\cdot A^{-1}\mathbf{u}\in L^1(\Om),\ \Div\mathbf{u}\in L^2(\Om)\right\}
\end{equation*}
(note that since $0<x_2<M$ on $\Om$, the condition $\mathbf{u}\cdot A^{-1}\mathbf{u}\in L^1(\Om)$ automatically implies $\mathbf{u}\in [L^2(\Om)]^2$). Its inner product and norm are defined by 
\begin{equation*}
	(\mathbf{u}, \mathbf{v})_{H_{\Div}(\Om;A^{-1})}=\int_\Om \left(\mathbf{u}\cdot A^{-1}\mathbf{v}+(\Div\mathbf{u})(\Div\mathbf{v})\right)\,\df x,\quad \|\mathbf{u}\|_{H_{\Div}(\Om;A^{-1})}=(\mathbf{u},\mathbf{u})_{H_{\Div}(\Om;A^{-1})}^\f{1}{2}. 
\end{equation*}
It is easily verified that the spaces 
\begin{equation*}
	H_{\Div}(\Om),\quad H_{\Div}(\Om;A^{-1})
\end{equation*}
are Hilbert spaces.

\begin{lemma}[{\cite[Theorem 3.18, p.~975]{Skrepek}, or \cite{Dautray}}]\label{11.22.L3}
	The space $[C_0^\iy(\R^2)]^2|_{\Om}$ is dense in $H_{\Div}(\Om)$. 
\end{lemma}

Now, we introduce the normal trace operator
\begin{equation}\label{11.22.6}
	\ga_\nu\colon H_{\Div}(\Om)\ra \left(H_{00}^\f{1}{2}(\Ga_2^0)\right)'
\end{equation}
as the unique continuous extension of
\begin{equation*}
	\langle \ga_\nu(\mathbf{u}), \phi \rangle_{(H_{00}^\f{1}{2}(\Ga_2^0))',\, H_{00}^\f{1}{2}(\Ga_2^0)} = -\int_{\Ga_2^0} u_2\phi\,\df x_1, \text{ for all } \mathbf{u}=(u_1,u_2)\in [C^\iy(\ol\Om)]^2 \text{ and } \phi\in H_{00}^\f{1}{2}(\Ga_2^0),
\end{equation*}
where $(H_{00}^\f{1}{2}(\Ga_2^0))'$ denotes the topological dual space of $H_{00}^\f{1}{2}(\Ga_2^0)$ with respect to the pivot space $L^2(\Ga_2^0)$.

\begin{lemma}\label{11.22.L4}
	The normal trace operator $\ga_\nu$ defined in \eqref{11.22.6} is a bounded linear operator. 
\end{lemma}

\begin{proof}
	First, note that the trace operator $\ga\colon H_\Ga^1(\Om)\ra H_{00}^\f{1}{2}(\Ga_2^0)$ defined in \eqref{11.22.5} is surjective and admits a bounded linear right inverse $E\colon H_{00}^\f{1}{2}(\Ga_2^0)\ra H_\Ga^1(\Om)$. 
	
	Now let $\mathbf{u}=(u_1,u_2)\in [C^\iy(\ol\Om)]^2$ and $\phi\in H_{00}^\f{1}{2}(\Ga_2^0)$ be arbitrary, and let $v=E\phi\in H_\Ga^1(\Om)$ be the lifting constructed above. Since the outward unit normal on $\Ga_2^0=(a_0,b_0)\ts\{0\}$ is $\nu=(0,-1)$ and $v=0$ on $\Ga=\pt\Om-\Ga_2^0$, integrating by parts yields
	\begin{equation*}
		\int_\Om (\Div\mathbf{u})v\,\df x + \int_\Om \mathbf{u}\cdot\nabla v\,\df x = \int_{\pt\Om} (\mathbf{u}\cdot\nu)v\,\df S = -\int_{\Ga_2^0} u_2\phi\,\df x_1.
	\end{equation*}
	Applying the Cauchy-Schwarz inequality, we deduce that
	\begin{equation*}
		\begin{split}
			\left|\int_{\Ga_2^0} u_2\phi\,\df x_1\right|
			&\leq \|\Div\mathbf{u}\|_{L^2(\Om)}\|v\|_{L^2(\Om)} + \|\mathbf{u}\|_{[L^2(\Om)]^2}\|\nabla v\|_{[L^2(\Om)]^2}\\
			&\leq \|\mathbf{u}\|_{H_{\Div}(\Om)}\|v\|_{H^1(\Om)} \leq C_1\|\mathbf{u}\|_{H_{\Div}(\Om)}\|\phi\|_{H_{00}^\f{1}{2}(\Ga_2^0)}.
		\end{split}
	\end{equation*}
	Taking the supremum over all $\phi\in H_{00}^\f{1}{2}(\Ga_2^0)$ with $\|\phi\|_{H_{00}^\f{1}{2}(\Ga_2^0)}\leq 1$, we obtain
	\begin{equation*}
		\|\ga_\nu(\mathbf{u})\|_{(H_{00}^\f{1}{2}(\Ga_2^0))'} \leq C_1\|\mathbf{u}\|_{H_{\Div}(\Om)},  \text{ for all } \mathbf{u}\in [C^\iy(\ol\Om)]^2.
	\end{equation*}
	By Lemma \ref{11.22.L3} (the density of $[C^\iy(\ol\Om)]^2$ in $H_{\Div}(\Om)$) and linearity, $\ga_\nu$ extends uniquely to a bounded linear operator from $H_{\Div}(\Om)$ into $(H_{00}^\f{1}{2}(\Ga_2^0))'$. We prove the lemma. 
\end{proof}

\subsection{Green's formulas}\label{S2.3}

Define a bilinear form $\mcT\colon H_{\Div}(\Om)\ts H_{\Ga}^1(\Om)\ra \R$ by
\begin{equation}\label{11.22.7}
	\mcT(\mathbf{u}, v)=\int_\Om \left((\Div\mathbf{u}) v+\mathbf{u}\cdot\nabla v\right)\,\df x,  \text{ for all } (\mathbf{u},v)\in H_{\Div}(\Om)\ts H_{\Ga}^1(\Om).
\end{equation}

\begin{lemma}[Green's formula]\label{11.22.L5}
	For all $\mathbf{u}\in H_{\Div}(\Om)$ and all $v\in H_\Ga^1(\Om)$, we have 
	\begin{equation}\label{11.22.8}
		\mcT(\mathbf{u}, v)=\llg \ga_\nu\mathbf{u}, \ga v\rrg_{(H_{00}^\f{1}{2}(\Ga_2^0))', H_{00}^\f{1}{2}(\Ga_2^0)}. 
	\end{equation}
\end{lemma}

\begin{proof}
	First, let $\mathbf{u}=(u_1,u_2)\in [C^\iy(\ol\Om)]^2$ and $v\in C_\Ga^\iy(\ol\Om)$. Since $\Om\s \R_+^2$ is a bounded convex  domain  with $C^3$ boundary $\pt\Om = \Ga \cup \Ga_2^0$, where $v=0$ on $\Ga$ and the outward unit normal on $\Ga_2^0=(a_0,b_0)\ts\{0\}$ is $\nu=(0,-1)$, the classical divergence theorem yields
	\begin{equation*}
		\begin{split}
			\mcT(\mathbf{u}, v)
			&=\int_\Om \left((\Div\mathbf{u}) v+\mathbf{u}\cdot\nabla v\right)\,\df x
			=\int_{\pt\Om} (\mathbf{u}\cdot\nu)v\,\df S\\
			&=-\int_{\Ga_2^0} u_2 v|_{\Ga_2^0}\,\df x_1
			=\llg \ga_\nu\mathbf{u}, \ga v\rrg_{(H_{00}^\f{1}{2}(\Ga_2^0))', H_{00}^\f{1}{2}(\Ga_2^0)},
		\end{split}
	\end{equation*}
	which proves \eqref{11.22.8} for all $(\mathbf{u},v)\in [C^\iy(\ol\Om)]^2 \ts C_\Ga^\iy(\ol\Om)$.
	
	Next, observe that both sides of \eqref{11.22.8} define continuous bilinear forms on $H_{\Div}(\Om)\ts H_\Ga^1(\Om)$. Indeed, by the Cauchy-Schwarz inequality,
	\begin{equation*} 
		|\mcT(\mathbf{u}, v)| \leq \|\Div\mathbf{u}\|_{L^2(\Om)}\|v\|_{L^2(\Om)} + \|\mathbf{u}\|_{[L^2(\Om)]^2}\|\nabla v\|_{[L^2(\Om)]^2} \leq \|\mathbf{u}\|_{H_{\Div}(\Om)}\|v\|_{H^1(\Om)},
	\end{equation*}
	while by Lemma \ref{11.22.L2} and Lemma \ref{11.22.L4},
	\begin{equation*}
		\left|\llg \ga_\nu\mathbf{u}, \ga v\rrg_{(H_{00}^\f{1}{2}(\Ga_2^0))', H_{00}^\f{1}{2}(\Ga_2^0)}\right|
		\leq \|\ga_\nu\mathbf{u}\|_{(H_{00}^\f{1}{2}(\Ga_2^0))'} \|\ga v\|_{H_{00}^\f{1}{2}(\Ga_2^0)}
		\leq C\|\mathbf{u}\|_{H_{\Div}(\Om)}\|v\|_{H^1(\Om)}.
	\end{equation*}
	Since $[C^\iy(\ol\Om)]^2$ is dense in $H_{\Div}(\Om)$ (by Lemma \ref{11.22.L3}) and $C_\Ga^\iy(\ol\Om)$ is dense in $H_\Ga^1(\Om)$ by definition, \eqref{11.22.8} holds for all $(\mathbf{u},v)\in H_{\Div}(\Om)\ts H_\Ga^1(\Om)$ by continuity and density. We have proved this lemma. 
\end{proof}

For each $\de\in (0,\de_0)$, we denote 
\begin{equation}\label{11.22.9}
	\wh\Om_\de=\Om\cap \{x_2>\de\}   \text{ and }   \Ga_2^\de=\Om\cap \{x_2=\de\},
\end{equation}
and set $\Ga_\de = \pt\wh\Om_\de - \Ga_2^\de = \pt\Om\cap\{x_2\geq \de\}$. Analogously to $H_{\Ga}^1(\Om)$ and $H_{00}^\f{1}{2}(\Ga_2^0)$, we denote by $H_{\Ga_\de}^1(\wh\Om_\de)$ the closure in $H^1(\wh\Om_\de)$ of $C_{\Ga_\de}^\iy(\ol{\wh\Om_\de})=\{u\in C^\iy(\ol{\wh\Om_\de})\colon \dist(\supp u, \Ga_\de)>0\}$, and define the trace operator
\begin{equation}\label{11.22.10}
	\ga^\de\colon H_{\Ga_\de}^1(\wh\Om_\de)\ra H_{00}^\f{1}{2}(\Ga_2^\de), \quad \ga^\de(u)=u|_{\Ga_2^\de} \text{ for all } u\in C_{\Ga_\de}^\iy(\ol{\wh\Om_\de}),
\end{equation}
and the normal trace operator 
\begin{equation}\label{11.22.11}
	\ga_\nu^\de\colon H_{\Div}(\wh\Om_\de)\ra \left(H_{00}^\f{1}{2}(\Ga_2^\de)\right)', \quad \ga_\nu^\de \mathbf{u}=-u_2|_{\Ga_2^\de}  \text{ for all } \mathbf{u}=(u_1,u_2)\in [C^\iy(\ol{\wh\Om_\de})]^2,
\end{equation}
extended by continuity. Additionally, we define the bilinear form
\begin{equation}\label{11.22.12}
	\mcT^\de\colon H_{\Div}(\wh\Om_\de)\ts H_{\Ga_\de}^1(\wh\Om_\de)\ra \R, \quad \mcT^\de(\mathbf{u}, v)=\int_{\wh\Om_\de} \left((\Div \mathbf{u})v+\mathbf{u}\cdot \nabla v\right)\,\df x. 
\end{equation}
By the exact same arguments as in Lemmas \ref{11.22.L2}, \ref{11.22.L4}, and \ref{11.22.L5}, we have:
\begin{enumerate}
	\item[{\rm (1)}] $\ga^\de$ is a bounded linear operator;
	\item[{\rm (2)}] $\ga_\nu^\de$ is a bounded linear operator;
	\item[{\rm (3)}] for all $\mathbf{u}\in H_{\Div}(\wh\Om_\de)$ and all $v\in H_{\Ga_\de}^1(\wh\Om_\de)$,
	\begin{equation*}
		\mcT^\de(\mathbf{u}, v)=\llg \ga_\nu^\de(\mathbf{u}), \ga^\de(v)\rrg_{(H_{00}^\f{1}{2}(\Ga_2^\de))', H_{00}^\f{1}{2}(\Ga_2^\de)}. 
	\end{equation*}
\end{enumerate}

\noindent Combining Lemma \ref{11.22.L5} with {\rm (3)} above, we obtain the following result.

\begin{corollary}\label{11.22.C1}
	For all $\mathbf{u}\in H_{\Div}(\Om)$ and all $v\in H_{\Ga}^1(\Om)$, we have 
	\begin{equation*}
		\llg \ga_\nu(\mathbf{u}), \ga (v)\rrg_{(H_{00}^\f{1}{2}(\Ga_2^0))', H_{00}^\f{1}{2}(\Ga_2^0)}=\lim_{\de\ra 0^+}\llg \ga_\nu^\de(\mathbf{u}|_{\wh\Om_\de}), \ga^\de(v|_{\wh\Om_\de})\rrg_{(H_{00}^\f{1}{2}(\Ga_2^\de))', H_{00}^\f{1}{2}(\Ga_2^\de)}. 
	\end{equation*}
\end{corollary}

\begin{proof}
	Let $\mathbf{u}\in H_{\Div}(\Om)$ and $v\in H_\Ga^1(\Om)$. Note that for every $\de\in (0,\de_0)$, we have $\mathbf{u}|_{\wh\Om_\de}\in H_{\Div}(\wh\Om_\de)$ and $v|_{\wh\Om_\de}\in H_{\Ga_\de}^1(\wh\Om_\de)$. By {\rm (3)}, the dominated convergence theorem, and Lemma \ref{11.22.L5}, we deduce that
	\begin{equation*}
		\begin{split}
			&\lim_{\de\ra 0^+}\llg \ga_\nu^\de(\mathbf{u}|_{\wh\Om_\de}), \ga^\de(v|_{\wh\Om_\de})\rrg_{(H_{00}^\f{1}{2}(\Ga_2^\de))', H_{00}^\f{1}{2}(\Ga_2^\de)}\\
			&=\lim_{\de\ra 0^+}\int_{\wh\Om_\de} \left((\Div \mathbf{u})v+\mathbf{u}\cdot \nabla v\right)\,\df x=\int_{\Om} \left((\Div \mathbf{u})v+\mathbf{u}\cdot \nabla v\right)\,\df x
			=\llg \ga_\nu(\mathbf{u}), \ga (v)\rrg_{(H_{00}^\f{1}{2}(\Ga_2^0))', H_{00}^\f{1}{2}(\Ga_2^0)}.
		\end{split}
	\end{equation*}
	This completes the proof.
\end{proof}

\begin{lemma}\label{11.22.L6}
	For each $\mathbf{u}\in H_{\Div}(\Om)$ and $v\in C_\Ga^\iy(\ol\Om)$, the functions
	\begin{equation}\label{11.22.14}
		F(v;\de)=\int_{\Ga_2^\de}|\ga^\de(v|_{\wh\Om_\de})|^2\,\df x_1, \quad  G(\mathbf{u}, v; \de)=\llg \ga_\nu^\de (\mathbf{u}|_{\wh\Om_\de}), \ga^\de (v|_{\wh\Om_\de})\rrg_{(H_{00}^\f{1}{2}(\Ga_2^\de))', H_{00}^\f{1}{2}(\Ga_2^\de)}
	\end{equation}
	are well-defined and continuous on $[0,\de_0]$ (with the convention $\wh\Om_0=\Om$, $\ga^0=\ga$, and $\ga_\nu^0=\ga_\nu$). In particular, both $F(v;\de)$ and $G(\mathbf{u}, v;\de)$ are continuous at $\de=0$. 
\end{lemma}

\begin{proof}
	Let $\mathbf{u}\in H_{\Div}(\Om)$ and $v\in C_\Ga^\iy(\ol\Om)$. Since $\dist(\supp v, \ol\Ga)>0$ and $\Ga_\de = \pt\Om\cap\{x_2\geq \de\}\s \Ga$, we have $\supp v \cap \Ga_\de = \emptyset$ for all $\de\in (0,\de_0]$. Hence, $v|_{\wh\Om_\de}\in H_{\Ga_\de}^1(\wh\Om_\de)$ for every $\de\in (0,\de_0]$, so both $F(v;\de)$ and $G(\mathbf{u},v;\de)$ are well-defined on $[0,\de_0]$.
	
	First, write $\Ga_2^\de = (a_\de, b_\de)\ts \{\de\}$ for $\de\in [0,\de_0]$. By the convexity and $C^3$-regularity of $\pt\Om$, the endpoints $a_\de$ and $b_\de$ depend continuously on $\de\in [0,\de_0]$. Since $v\in C^\iy(\ol\Om)$ is bounded and continuous on $\ol\Om$, we have
	\begin{equation*}
		F(v;\de) = \int_{a_\de}^{b_\de} |v(x_1,\de)|^2\,\df x_1, \quad \de\in [0,\de_0],
	\end{equation*}
	which is continuous on $[0,\de_0]$ by the continuity of $v$, $a_\de$, and $b_\de$.
	
	Second, by Lemma \ref{11.22.L5} (for $\de=0$) and property {\rm (3)} following \eqref{11.22.12} (for $\de\in (0,\de_0]$), we obtain
	\begin{equation*}
		G(\mathbf{u},v;\de) = \int_{\wh\Om_\de} \left((\Div\mathbf{u})v + \mathbf{u}\cdot\nabla v\right)\,\df x,  \text{ for all } \de\in [0,\de_0].
	\end{equation*}
	Since $(\Div\mathbf{u})v + \mathbf{u}\cdot\nabla v \in L^1(\Om)$, the absolute continuity of the Lebesgue integral implies that $G(\mathbf{u},v;\de)$ is continuous on $[0,\de_0]$. We have proved the lemma. 
\end{proof}

To establish Green's formula for the weighted space $H_{\Div}(\Om;A^{-1})$, we introduce the bilinear form
\begin{equation}\label{11.22.15}
	\mcT_A\colon H_{\Div}(\Om;A^{-1})\ts H_\Ga^1(\Om;A)\ra \R 
\end{equation}
defined by 
\begin{equation*}
	\mcT_A(\mathbf{u},v)=\int_\Om \left((\Div\mathbf{u})v+\mathbf{u}\cdot \nabla v\right)\,\df x, \text{ for all } \mathbf{u}\in H_{\Div}(\Om;A^{-1}) \text{ and } v\in H_\Ga^1(\Om;A). 
\end{equation*}

\begin{lemma}\label{11.22.L7}
	Let $\al\in [1,2)$. For any $\mathbf{u}\in H_{\Div}(\Om; A^{-1})$, we have 
	\begin{equation*}
		\mcT_A(\mathbf{u},v)=0, \text{ for all } v\in H_\Ga^1(\Om;A). 
	\end{equation*}
\end{lemma}

\begin{proof}
	First, note that $\mcT_A$ is a continuous bilinear form on $H_{\Div}(\Om; A^{-1})\ts H_\Ga^1(\Om;A)$. Indeed, writing $\mathbf{u}\cdot\nabla v = (A^{-\f{1}{2}}\mathbf{u})\cdot (A^{\f{1}{2}}\nabla v)$ and applying the Cauchy-Schwarz inequality, we obtain
	\begin{equation*}
		\begin{split} 
			|\mcT_A(\mathbf{u},v)|
			&\leq \|\Div\mathbf{u}\|_{L^2(\Om)}\|v\|_{L^2(\Om)} + \|A^{-\f{1}{2}}\mathbf{u}\|_{[L^2(\Om)]^2}\|A^{\f{1}{2}}\nabla v\|_{[L^2(\Om)]^2}\\
			&\leq C\|\mathbf{u}\|_{H_{\Div}(\Om;A^{-1})}\|v\|_{H_\Ga^1(\Om;A)},
		\end{split}
	\end{equation*}
	where we used the Poincar\'e inequality \eqref{11.22.4}. Hence, by the density of $C_\Ga^\iy(\ol\Om)$ in $H_\Ga^1(\Om;A)$, it suffices to show that $\mcT_A(\mathbf{u},v)=0$ for all $v\in C_\Ga^\iy(\ol\Om)$.
	
	Fix $\mathbf{u}=(u_1,u_2)\in H_{\Div}(\Om;A^{-1})$ and $v\in C_\Ga^\iy(\ol\Om)$. Since $H_{\Div}(\Om;A^{-1})\hra H_{\Div}(\Om)$ and $v\in C_\Ga^\iy(\ol\Om)\s H_\Ga^1(\Om)$, it follows from Lemma \ref{11.22.L5} and Corollary \ref{11.22.C1} that
	\begin{equation}\label{11.22.16}
		\mcT_A(\mathbf{u},v)=\mcT(\mathbf{u},v)=\lim_{\de\ra 0^+}\llg \ga_\nu^\de(\mathbf{u}|_{\wh\Om_\de}), \ga^\de(v|_{\wh\Om_\de})\rrg_{(H_{00}^\f{1}{2}(\Ga_2^\de))', H_{00}^\f{1}{2}(\Ga_2^\de)}.
	\end{equation}
	On the other hand, since $\mathbf{u}\cdot A^{-1}\mathbf{u} = u_1^2 + x_2^{-\al}u_2^2 \in L^1(\Om)$, Fubini's theorem implies that for almost every $s\in (0,\de_0)$, we have $u_2(\cdot, s)\in L^2(\Ga_2^s)$ and therefore
	\begin{equation*}
		\llg \ga_\nu^s(\mathbf{u}|_{\wh\Om_s}), \ga^s(v|_{\wh\Om_s})\rrg_{(H_{00}^\f{1}{2}(\Ga_2^s))', H_{00}^\f{1}{2}(\Ga_2^s)} = -\int_{\Ga_2^s} u_2(x_1,s)v(x_1,s)\,\df x_1.
	\end{equation*}
	(Indeed, this holds for $\mathbf{u}\in [C^\iy(\ol\Om)]^2$, and extends to $\mathbf{u}\in H_{\Div}(\Om;A^{-1})$ for a.e. $s\in (0,\de_0)$ by approximating $\mathbf{u}$ in $H_{\Div}(\Om)$ via Lemma \ref{11.22.L3} and using property {\rm (3)} following \eqref{11.22.12}.)
	
	Now, for any $\de\in (0,\de_0)$, integrating over $s\in (0,\de)$ and applying the Cauchy-Schwarz inequality yields
	\begin{equation*}
		\begin{split}
			&\frac{1}{\de}\int_0^\de \left|\llg \ga_\nu^s(\mathbf{u}|_{\wh\Om_s}), \ga^s(v|_{\wh\Om_s})\rrg_{(H_{00}^\f{1}{2}(\Ga_2^s))', H_{00}^\f{1}{2}(\Ga_2^s)}\right|\df s\\
			&\leq \frac{1}{\de}\int_0^\de \int_{\Ga_2^s} |u_2(x_1,s)|\,|v(x_1,s)|\,\df x_1\,\df s\leq \frac{1}{\de}\left(\int_0^\de \int_{\Ga_2^s} s^{-\al}u_2^2\,\df x_1\,\df s\right)^{\f{1}{2}} \left(\int_0^\de s^\al F(v;s)\,\df s\right)^{\f{1}{2}},
		\end{split}
	\end{equation*}
	where $F(v;s)=\int_{\Ga_2^s}|\ga^s(v|_{\wh\Om_s})|^2\,\df x_1$. By Lemma \ref{11.22.L6}, $F(v;\cdot)$ is continuous on $[0,\de_0]$, so $M_v:=\sup_{s\in [0,\de_0]} F(v;s)<+\iy$. Hence,
	\begin{equation*}
		\frac{1}{\de}\int_0^\de \left|\llg \ga_\nu^s(\mathbf{u}|_{\wh\Om_s}), \ga^s(v|_{\wh\Om_s})\rrg_{(H_{00}^\f{1}{2}(\Ga_2^s))', H_{00}^\f{1}{2}(\Ga_2^s)}\right|\df s
		\leq \left(\frac{M_v}{\al+1}\right)^{\f{1}{2}} \de^{\f{\al-1}{2}} \left(\int_{\Om-\wh\Om_\de} x_2^{-\al}u_2^2\,\df x\right)^{\f{1}{2}}.
	\end{equation*}
	Since $\al\geq 1$, we have $\de^{\f{\al-1}{2}}\leq 1$ for $\de\in (0,1)$, while $\int_{\Om-\wh\Om_\de} x_2^{-\al}u_2^2\,\df x \to 0$ as $\de\to 0^+$ by the dominated convergence theorem. Consequently, the right-hand side converges to $0$ as $\de\to 0^+$. By the mean value theorem for integrals (or selecting a sequence $\de_k\to 0^+$), there exists a sequence $\de_k\to 0^+$ such that
	\begin{equation*}
		\lim_{k\to\iy}\llg \ga_\nu^{\de_k}(\mathbf{u}|_{\wh\Om_{\de_k}}), \ga^{\de_k}(v|_{\wh\Om_{\de_k}})\rrg_{(H_{00}^\f{1}{2}(\Ga_2^{\de_k}))', H_{00}^\f{1}{2}(\Ga_2^{\de_k})} = 0.
	\end{equation*}
	Combining this with \eqref{11.22.16} yields $\mcT_A(\mathbf{u},v)=0$.
\end{proof}

\begin{lemma}\label{11.22.L8}
	The space $[C^\iy(\ol\Om)]^2\cap H_{\Div}(\Om;A^{-1})$ is dense in $H_{\Div}(\Om;A^{-1})$. 
\end{lemma}

\begin{proof}
	Set $\mcV := [C^\iy(\ol\Om)]^2\cap H_{\Div}(\Om;A^{-1})$. Since $H_{\Div}(\Om;A^{-1})$ is a Hilbert space, by the Hahn-Banach (or Riesz representation) theorem, it suffices to show that any $\mathbf{u}\in H_{\Div}(\Om;A^{-1})$ satisfying
	\begin{equation}\label{11.22.17}
		(\mathbf{u},\mathbf{w})_{H_{\Div}(\Om;A^{-1})} = \int_\Om \left(\mathbf{u}\cdot A^{-1}\mathbf{w} + (\Div\mathbf{u})(\Div\mathbf{w})\right)\df x = 0, \text{ for all } \mathbf{w}\in \mcV
	\end{equation}
	must vanish identically, i.e., $\mathbf{u}=\mathbf{0}$.
	
	Let $\mathbf{u}\in H_{\Div}(\Om;A^{-1})$ satisfy \eqref{11.22.17}. First, taking $\mathbf{w}\in [C_0^\iy(\Om)]^2\s \mcV$ in \eqref{11.22.17}, we obtain in the sense of distributions that
	\begin{equation*}
		\nabla(\Div\mathbf{u}) = A^{-1}\mathbf{u}  \text{ in } [(C_0^\iy(\Om))']^2.
	\end{equation*}
	Set $v := \Div\mathbf{u}\in L^2(\Om)$. Then $A\nabla v = \mathbf{u}$ in $\Om$, and hence
	\begin{equation*}
		\int_\Om \nabla v\cdot A\nabla v\,\df x = \int_\Om \mathbf{u}\cdot A^{-1}\mathbf{u}\,\df x < +\iy,
	\end{equation*}
	which shows that $v\in \wh H^1(\Om;A)=H^1(\Om;A)$ by Lemma \ref{11.22.L9}.

	For any $\mathbf{w}\in [C^\iy(\ol\Om)]^2$ with $\dist(\supp\mathbf{w}, \ol{\Ga_2^0})>0$, the weight $x_2^{-\al}$ is bounded on $\supp\mathbf{w}$, so $\mathbf{w}\in \mcV$. Substituting such $\mathbf{w}$ into \eqref{11.22.17} and using $\mathbf{u}=A\nabla v$ and $\Div\mathbf{u}=v$, we obtain
	\begin{equation*}
		0 = \int_\Om \left(\nabla v\cdot \mathbf{w} + v\,\Div\mathbf{w}\right)\df x.
	\end{equation*}
	Since $\Ga = \pt\Om-\Ga_2^0$ is of class $C^3$ (away from $\ol{\Ga_2^0}$) and $v\in H_{\mathrm{loc}}^1(\ol\Om-\ol{\Ga_2^0})$, integrating by parts yields $\int_\Ga v(\mathbf{w}\cdot\nu)\,\df S = 0$ for all such $\mathbf{w}$, which implies $v=0$ on $\Ga$, and hence $v\in H_\Ga^1(\Om;A)$.
	
	Finally, applying Lemma \ref{11.22.L7} to $\mathbf{u}\in H_{\Div}(\Om;A^{-1})$ and $v\in H_\Ga^1(\Om;A)$, we conclude that
	\begin{equation*}
		0 = \mcT_A(\mathbf{u},v) = \int_\Om \left((\Div\mathbf{u})v + \mathbf{u}\cdot\nabla v\right)\df x = \int_\Om \left(|\Div\mathbf{u}|^2 + \mathbf{u}\cdot A^{-1}\mathbf{u}\right)\df x = \|\mathbf{u}\|_{H_{\Div}(\Om;A^{-1})}^2.
	\end{equation*}
	Thus $\mathbf{u}=\mathbf{0}$, which completes the proof.
\end{proof}

\begin{theorem}[Green's formula]\label{11.23.T1}
	Let $\al\in [1,2)$. Then the following assertions hold:
	
	{\rm (i)} $H_{\Div}(\Om;A^{-1})\s \ker \ga_\nu=\left\{\mathbf{u}\in H_{\Div}(\Om)\colon \ga_\nu(\mathbf{u})=0\right\}$. 
	
	{\rm (ii)} For every $u\in D(\mcA)$, we have $A\nabla u\in H_{\Div}(\Om;A^{-1})$, and the conormal derivative of $u$ on $\Ga_2^0$ vanishes in the sense that
	\begin{equation*}
		\f{\pt u}{\pt \nu_A}:=\ga_\nu(A\nabla u)=0 \quad \text{in }\left(H_{00}^\f{1}{2}(\Ga_2^0)\right)'.
	\end{equation*}
	Furthermore, for all $(u,v)\in D(\mcA)\ts H_\Ga^1(\Om;A)$, the following Green's formula holds: 
	\begin{equation}\label{11.23.Green}
		\int_\Om (\mcA u)v\,\df x = -\int_\Om [\Div(A\nabla u)]v\,\df x=\int_\Om \nabla u\cdot A\nabla v\,\df x. 
	\end{equation}
\end{theorem}

\begin{proof}
	{\rm (i)} Let $\mathbf{u}\in H_{\Div}(\Om;A^{-1})$. Since $0<x_2<M$ on $\Om$, we have $H_{\Div}(\Om;A^{-1})\s H_{\Div}(\Om)$. For any $\phi\in H_{00}^\f{1}{2}(\Ga_2^0)$, as shown in the proof of Lemma \ref{11.22.L4}, there exists $v\in H_\Ga^1(\Om)$ such that $\ga v = \phi$. Since $H_\Ga^1(\Om)\s H_\Ga^1(\Om;A)$, applying Lemma \ref{11.22.L5} and Lemma \ref{11.22.L7} gives
	\begin{equation*}
		\llg \ga_\nu(\mathbf{u}), \phi\rrg_{(H_{00}^\f{1}{2}(\Ga_2^0))', H_{00}^\f{1}{2}(\Ga_2^0)}
		= \llg \ga_\nu(\mathbf{u}), \ga v\rrg_{(H_{00}^\f{1}{2}(\Ga_2^0))', H_{00}^\f{1}{2}(\Ga_2^0)}
		= \mcT(\mathbf{u}, v) = \mcT_A(\mathbf{u}, v) = 0.
	\end{equation*}
	Since $\phi\in H_{00}^\f{1}{2}(\Ga_2^0)$ is arbitrary, we conclude that $\ga_\nu(\mathbf{u})=0$, i.e., $\mathbf{u}\in \ker\ga_\nu$.
	
	{\rm (ii)} Let $u\in D(\mcA)$. By the definition of $D(\mcA)$, we have $u\in H_\Ga^1(\Om;A)$ and $\Div(A\nabla u) = -\mcA u \in L^2(\Om)$. Setting $\mathbf{u} := A\nabla u$, we observe that $\mathbf{u}\cdot A^{-1}\mathbf{u} = \nabla u\cdot A\nabla u \in L^1(\Om)$ and $\Div\mathbf{u} \in L^2(\Om)$, which implies $\mathbf{u} = A\nabla u \in H_{\Div}(\Om;A^{-1})$ (since $0<x_2<M$ on $\Om$). By part {\rm (i)}, it immediately follows that $\f{\pt u}{\pt\nu_A} = \ga_\nu(A\nabla u) = 0$ in $(H_{00}^\f{1}{2}(\Ga_2^0))'$.
	
	Finally, for any $u\in D(\mcA)$ and $v\in H_\Ga^1(\Om;A)$, applying Lemma \ref{11.22.L7} with $\mathbf{u} = A\nabla u \in H_{\Div}(\Om;A^{-1})$ yields
	\begin{equation*}
		0 = \mcT_A(A\nabla u, v) = \int_\Om \left([\Div(A\nabla u)]v + A\nabla u\cdot\nabla v\right)\df x,
	\end{equation*}
	which rearranges to \eqref{11.23.Green} because $A$ is symmetric. This completes the proof.
\end{proof}

\section{Shape design}\label{S3} 

In this section, we introduce a shape design method for mixed Dirichlet-Neumann boundary value problems.

\subsection{Preliminary}\label{S3.1}

Let $\de_0>0$ be given in \eqref{11.22.3}. For each $\de\in (0,\de_0 )$, we choose a bounded convex domain $\Om_\de\s \Om$ with $\pt\Om_\de\in C^3$ such that
\begin{equation}\label{11.23.1}
	\{x\in\Om\colon x_2>2\de\}\s \Om_\de\s \{x\in\Om\colon x_2>\de\}, \quad \Om_{\de_1}\s \Om_{\de_2} \mbox{ for all } 0\leq \de_2\leq \de_1\leq \de_0. 
\end{equation}
Denote 
\begin{equation*}
	\Ga_\de=\Ga\cap \pt\Om_\de \text{ and }  \Ga_2^\de=\pt\Om_\de\cap \Om
\end{equation*}
(note that $\Ga_2^\de\neq \Om\cap\{x_2=\de\}$ here, as $\pt\Om_\de$ is of class $C^3$, and $\pt\Om_\de = \Ga_\de \cup \Ga_2^\de$ with $\Ga_\de\cap\Ga_2^\de=\emptyset$). Set 
\begin{equation*}
	A_\de=A|_{\Om_\de}, \quad \mcA_\de u=-\Div(A_\de\nabla u)  \text{ in }\Om_\de, \quad Q_\de=\Om_\de \ts (0,T), \quad \Sigma_\de=\pt\Om_\de \ts (0,T). 
\end{equation*}

Consider the following uniformly parabolic approximating equation:
\begin{equation}\label{11.23.2}
	\begin{cases}
		\pt_t\vp_\de-\Div(A_\de\nabla \vp_\de)=f_\de, &\text{in }Q_\de, \\
		\vp_\de=0, &\text{on }\Ga_\de\ts (0,T),\\
		\f{\pt\vp_\de}{\pt\nu_{A_\de}}=0, &\text{on }\Ga_2^\de\ts (0,T),\\
		\vp_\de(0)=\vp_\de^0, &\text{in }\Om_\de,  
	\end{cases}
\end{equation}
where $\vp_\de^0\in L^2(\Om_\de)$, $f_\de\in L^2(Q_\de)$, and $\f{\pt\vp_\de}{\pt\nu_{A_\de}} = A_\de\nabla\vp_\de\cdot\nu_\de$ denotes the classical conormal derivative on $\Ga_2^\de$ (with $\nu_\de$ being the outward unit normal to $\pt\Om_\de$).

Analogously to Section \ref{S2.1}, we introduce the spaces
\begin{equation*}
	C_{\Ga_\de}^\iy(\ol\Om_\de), \quad H_{\Ga_\de}^1(\Om_\de;A_\de), \quad H_{\Ga_\de}^{-1}(\Om_\de;A_\de),\quad H^1(\Om_\de;A_\de), \quad \wh H^1(\Om_\de;A_\de), \quad H^2(\Om_\de;A_\de),\quad D(\mcA_\de);
\end{equation*}
specifically, 
\begin{equation*}
	\begin{split}
		C_{\Ga_\de}^\iy(\ol\Om_\de)
		&=\left\{u\in C^\iy(\ol\Om_\de)\colon \dist(\supp u, \ol{\Ga_\de})>0\right\},\\
		\wh H^1(\Om_\de;A_\de)
		&=\left\{u\in L^2(\Om_\de)\colon \int_{\Om_\de} \nabla u\cdot A_\de\nabla u\,\df x<+\iy\right\},\\
		D(\mcA_\de)
		&=\left\{u\in H_{\Ga_\de}^1(\Om_\de;A_\de)\colon \mcA_\de u\in L^2(\Om_\de), \ \f{\pt u}{\pt \nu_{A_\de}}=0 \text{ on } \Ga_2^\de\right\}, 
	\end{split}
\end{equation*}
and $H^1(\Om_\de;A_\de)$ is the closure of $C^\iy(\ol\Om_\de)$ in $\wh H^1(\Om_\de;A_\de)$. For each $\de\in (0,\de_0)$, since $\de^\al < x_2^\al < M^\al$ on $\Om_\de$ and $\pt\Om_\de\in C^3$, we have $\wh H^1(\Om_\de;A_\de)=H^1(\Om_\de;A_\de) = H^1(\Om_\de)$ and $H_{\Ga_\de}^1(\Om_\de;A_\de) = H_{\Ga_\de}^1(\Om_\de)$ with equivalent norms.

\begin{remark}\label{11.23.R1}
	Lemmas \ref{01.28.L3} and \ref{11.22.L1} as well as Remark \ref{11.22.R1} remain valid on $H_{\Ga_\de}^1(\Om_\de;A_\de)$ for every $\de\in (0,\de_0 )$, with all constants independent of $\de$. Indeed, by the convexity of $\Om_\de$, its bottom boundary $\Ga_2^\de$ can be represented as the graph of a convex function (extended upward by $\Ga_\de$), so any $u\in C_{\Ga_\de}^\iy(\ol\Om_\de)$ can be extended by zero across $\Ga_\de$ to the region above this graph while vanishing near $\{x_2=M\}$, and the outward unit normal $\nu_\de=(\nu_{\de,1},\nu_{\de,2})$ on $\Ga_2^\de$ satisfies $\nu_{\de,2}\leq 0$. Consequently, in the proof of Lemma \ref{01.28.L3} for $\al=1$, the divergence theorem gives
	\begin{equation*}
		\int_{\Om_\de} \frac{\pt}{\pt x_2}(x_2 u^2)\,\df x = \int_{\Ga_2^\de} x_2 u^2 \nu_{\de,2}\,\df S \leq 0,
	\end{equation*}
	and the same sign condition $\nu_{\de,2}\leq 0$ ensures that all boundary terms arising from integration by parts on $\Om_\de\cap\{x_2>\e\}$ ($\e>\de$) remain nonpositive. Thus, the exact same arguments apply without altering any constants.
\end{remark}

\begin{remark}\label{11.23.R2}
	From now on, we allow $\de\in [0,\de_0]$, with $\de=0$ corresponding to the original problem, i.e., $\Om_0=\Om$, $\Ga_0=\Ga$, $A_0=A$, $\mcA_0=\mcA$, and so on. 
\end{remark}

\begin{remark}
	In the degenerate case $\de=0$ (where $\mcA_0=\mcA$ and $\Om_0=\Om$), it follows from Theorem \ref{11.23.T1} that any function $u\in D(\mcA)$ implicitly satisfies the homogeneous conormal boundary condition $\f{\pt u}{\pt \nu_A}=0$ on the degenerate boundary $\Ga_2^0$ in the sense of traces (i.e., the integration-by-parts identity $\int_\Om (\mcA u)v\,\df x = \int_\Om \nabla u\cdot A\nabla v\,\df x$ holds for all $v\in H_\Ga^1(\Om;A)$).
\end{remark}

\subsection{Spectrum}\label{S3.2}

For each $\de\in [0,\de_0]$, by the Lax-Milgram theorem, Green's formula (Theorem \ref{11.23.T1}), the compact embedding $H_{\Ga_\de}^1(\Om_\de;A_\de)\hra L^2(\Om_\de)$ (Lemma \ref{11.22.L1}), and the Poincar\'e inequality \eqref{11.22.4} (Remark \ref{11.22.R1}), the operator $\mcA_\de\colon D(\mcA_\de)\s L^2(\Om_\de)\ra L^2(\Om_\de)$ is a positive self-adjoint operator with compact resolvent. Consequently, $\mcA_\de$ admits a sequence of positive eigenvalues (repeated according to their finite multiplicities)
\begin{equation}\label{11.23.3}
	0<\la_1^\de\leq \la_2^\de\leq \cdots\leq \la_n^\de\leq \cdots, \quad \la_n^\de \ra+\iy \text{ as } n\ra\iy. 
\end{equation}
Moreover, letting $C_P>0$ denote the $\de$-independent Poincar\'e constant in \eqref{11.22.4}, we have 
\begin{equation*}
	C_P^{-1}\leq \la_1^\de \text{ for all } \de\in [0,\de_0]. 
\end{equation*}
We denote by $\Phi_n^\de\in D(\mcA_\de)$ the $L^2(\Om_\de)$-normalized eigenfunction of $\mcA_\de$ associated with $\la_n^\de$. By the spectral theorem for compact self-adjoint operators, $\{\Phi_n^\de\}_{n\in\N^*}$ forms an orthonormal basis of $L^2(\Om_\de)$, and $\{(\la_n^\de)^{-\f{1}{2}}\Phi_n^\de\}_{n\in\N^*}$ forms an orthonormal basis of $H_{\Ga_\de}^1(\Om_\de;A_\de)$ (equipped with the inner product $(u,v)_{H_{\Ga_\de}^1(\Om_\de;A_\de)} = \int_{\Om_\de} \nabla u\cdot A_\de\nabla v\,\df x$).

\subsection{Weak solutions}\label{S3.3}

\begin{definition}\label{11.23.D1}
	Let $\de\in [0,\de_0]$ (see Remark \ref{11.23.R2}). A function 
	\begin{equation*}
		\vp_\de\in C([0,T]; L^2(\Om_\de))\cap L^2(0,T; H_{\Ga_\de}^1(\Om_\de;A_\de)),  \text{ with }  \pt_t\vp_\de \in L^2(0,T; H_{\Ga_\de}^{-1}(\Om_\de;A_\de))
	\end{equation*}
	is called a weak solution of \eqref{11.23.2} associated with $(\vp_\de^0, f_\de)$ if:
	\begin{enumerate}
		\item[{\rm (i)}] for each $v\in H_{\Ga_\de}^1(\Om_\de;A_\de)$ and a.e.~$t\in (0,T)$, 
		\begin{equation*}
			\llg \pt_t\vp_\de(t), v\rrg_{H_{\Ga_\de}^{-1}(\Om_\de;A_\de), H_{\Ga_\de}^1(\Om_\de;A_\de)}+(\vp_\de(t), v)_{H_{\Ga_\de}^1(\Om_\de;A_\de)}=(f_\de(t), v)_{L^2(\Om_\de)};
		\end{equation*}
		\item[{\rm (ii)}] $\vp_\de(0)=\vp_\de^0$ in $L^2(\Om_\de)$. 
	\end{enumerate}
\end{definition}

\begin{theorem}\label{11.23.T2}
	Let $\de\in [0,\de_0]$ and $\al\in [1,2)$. For each $\vp_\de^0\in L^2(\Om_\de)$ and $f_\de\in L^2(Q_\de)$, there exists a unique weak solution $\vp_\de$ to \eqref{11.23.2} associated with $(\vp_\de^0,f_\de)$ satisfying 
	\begin{equation}\label{11.23.4}
		\begin{split}
			&\|\vp_\de\|_{L^\iy(0,T; L^2(\Om_\de))}+\|\vp_\de\|_{L^2(0,T; H_{\Ga_\de}^1(\Om_\de;A_\de))}+\|\pt_t\vp_\de\|_{L^2(0,T; H_{\Ga_\de}^{-1}(\Om_\de;A_\de))}\\
			&\leq C\left(\|\vp_\de^0\|_{L^2(\Om_\de)}+\|f_\de\|_{L^2(Q_\de)}\right), 
		\end{split}
	\end{equation}
	where the positive constant $C$ depends only on $\al$, $M$, and $T$. 
	
	Moreover, if we further assume that $\vp_\de^0\in H_{\Ga_\de}^1(\Om_\de;A_\de)$ and $f_\de\in H^1(0,T; L^2(\Om_\de))$, then 
	\begin{equation}\label{11.23.5}
		\begin{split}
			&\|\vp_\de\|_{L^\iy(0,T; H_{\Ga_\de}^1(\Om_\de;A_\de))}+\|\vp_\de\|_{L^2(0,T; D(\mcA_\de))}+\|\vp_\de\|_{L^2(0,T; H^2(\Om_{4\de_0}))}+\|\pt_t\vp_\de\|_{L^2(Q_\de)}\\
			&\leq C\left(\|\vp_\de^0\|_{H_{\Ga_\de}^1(\Om_\de;A_\de)}+\|f_\de\|_{H^1(0,T; L^2(\Om_\de))}\right), 
		\end{split}
	\end{equation}
where $\Omega_{4\delta_0}\subset\Omega_\delta$ is defined as in \eqref{11.23.1} (with $\delta$ replaced by $4\delta_0$), and the positive constant $C$ depends only on $\alpha$, $M$, the fixed domain $D=\Omega_{\delta_0}$, and $T$, independently of $\delta\in[0,\delta_0]$.
\end{theorem}

We postpone the proof of this theorem to Section \ref{S5}.

\begin{proposition}\label{11.23.P1}
	Let $\de\in [0,\de_0]$, $\vp_\de^0\in L^2(\Om_\de)$, and $f_\de\in L^2(Q_\de)$. Then a function 
	\begin{equation*}
		\vp_\de\in C([0,T]; L^2(\Om_\de))\cap L^2(0,T; H_{\Ga_\de}^1(\Om_\de;A_\de)), \text{ with }  \pt_t\vp_\de\in L^2(0,T; H_{\Ga_\de}^{-1}(\Om_\de;A_\de))
	\end{equation*}
	is a weak solution to \eqref{11.23.2} associated with $(\vp_\de^0, f_\de)$ if and only if $\vp_\de\in C([0,T]; L^2(\Om_\de))\cap L^2(0,T; H_{\Ga_\de}^1(\Om_\de;A_\de))$ satisfies
	\begin{equation}\label{11.23.26}
		-\iint_{Q_\de} \vp_\de \pt_t\psi\,\df x\,\df t + \iint_{Q_\de} \nabla\vp_\de \cdot A_\de\nabla \psi\,\df x\,\df t = \iint_{Q_\de} f_\de \psi\,\df x\,\df t + \int_{\Om_\de}\vp_\de^0 \psi(\cdot,0)\,\df x
	\end{equation}
	for all $\psi\in C^\iy(\ol Q_\de)$ with $\dist(\supp \psi(t),\ol{\Ga_\de})>0$ for all $t\in [0,T]$ and $\psi(\cdot,T)=0$.  
\end{proposition}

We postpone the proof of this proposition to Section \ref{S5}.

\subsection{Shape design}\label{S3.4}

\begin{theorem}\label{11.23.T3}
	Let $\al\in [1,2)$, $\vp^0\in C_\Ga^\iy(\ol\Om)$, and $f\in H^1(0,T; L^2(\Om))$. For each $\de\in [0,\de_0]$, let $\vp_\de$ be the weak solution to \eqref{11.23.2} associated with $(\vp_\de^0, f_\de):=(\vp^0|_{\Om_\de}, f|_{Q_\de})$. Then, for each $\be\in (0,\de_0]$, as $\de\to 0^+$, 
	\begin{equation}\label{11.23.19}
		\begin{split}
			\vp_\de 
			&\ra \vp_0 \text{ weakly in } L^2(0,T; H^2(\Om_\be;A_\be)\cap H_{\Ga_\be}^1(\Om_\be;A_\be)), \\
			\vp_\de 
			&\ra \vp_0 \text{ strongly in } L^2(Q_\be).
		\end{split}
	\end{equation}
	Moreover, if we additionally assume that $f=0$, then, as $\de\to 0^+$, 
	\begin{equation}\label{11.23.20}
		\f{\pt\vp_{\de}}{\pt\nu}\ra \f{\pt\vp_0}{\pt\nu} \text{ strongly in } L^2(0,T; L^2(\Ga_{4\de_0})). 
	\end{equation}
\end{theorem}

\begin{proof}[Proof of Theorem \ref{11.23.T3}]
	We divide the proof into two steps.

	\textit{Step 1}. Proof of \eqref{11.23.19}. 
	
	Since $\vp^0\in C_\Ga^\iy(\ol\Om)$ and $\ol{\Ga_\de}\s \ol\Ga$, we have $\vp_\de^0=\vp^0|_{\Om_\de}\in C_{\Ga_\de}^\iy(\ol{\Om_\de})\s H_{\Ga_\de}^1(\Om_\de;A_\de)$ for all $\de\in [0,\de_0]$, with $\|\vp_\de^0\|_{H_{\Ga_\de}^1(\Om_\de;A_\de)}\leq \|\vp^0\|_{H_\Ga^1(\Om;A)}$. Moreover, $f_\de=f|_{Q_\de}\in H^1(0,T; L^2(\Om_\de))$ satisfies $\|f_\de\|_{H^1(0,T; L^2(\Om_\de))}\leq \|f\|_{H^1(0,T; L^2(\Om))}$. By Theorem \ref{11.23.T2} (estimates \eqref{11.23.4} and \eqref{11.23.5}), there exists a constant $C>0$, independent of $\de\in [0,\de_0]$, such that
	\begin{equation}\label{11.23.21}
		\begin{split}
			&\|\vp_\de\|_{L^\iy(0,T; H_{\Ga_\de}^1(\Om_\de;A_\de))} + \|\vp_\de\|_{L^2(0,T; D(\mcA_\de))} + \|\pt_t\vp_\de\|_{L^2(Q_\de)} \\
			&\leq C\left(\|\vp^0\|_{H_\Ga^1(\Om;A)} + \|f\|_{H^1(0,T; L^2(\Om))}\right).
		\end{split}
	\end{equation}

	Fix any $\be\in (0,\de_0]$. For every $\de\in (0,\be)$, we have $\Om_\be\s \Om_\de$, $\Ga_\be\s \Ga_\de$, $A_\be=A_\de|_{\Om_\be}$, and $\mcA_\be(\vp_\de|_{\Om_\be})=(\mcA_\de\vp_\de)|_{\Om_\be}\in L^2(Q_\be)$. By the definition of $H^2(\Om_\be;A_\be)\cap H_{\Ga_\be}^1(\Om_\be;A_\be)$ together with \eqref{11.23.21}, the family $\{\vp_\de\}_{0<\de<\be}$ is uniformly bounded in
	\begin{equation*}
		L^2(0,T; H^2(\Om_\be;A_\be)\cap H_{\Ga_\be}^1(\Om_\be;A_\be)) \cap H^1(0,T; L^2(\Om_\be)).
	\end{equation*}
	On the fixed subdomain $\Om_\be$, the diffusion matrix $A_\be$ is uniformly elliptic and bounded, so $H^2(\Om_\be;A_\be)\cap H_{\Ga_\be}^1(\Om_\be;A_\be)$ embeds continuously into $H^1(\Om_\be)$. By the compact embedding $H^1(\Om_\be)\hookrightarrow\hookrightarrow L^2(\Om_\be)$ and the Aubin-Lions lemma, for any sequence $\de_k\to 0^+$, there exist a subsequence (still denoted by $\{\de_k\}$) and a function
	\begin{equation*}
		\wt\vp_\be \in L^2(0,T; H^2(\Om_\be;A_\be)\cap H_{\Ga_\be}^1(\Om_\be;A_\be)) \cap H^1(0,T; L^2(\Om_\be))
	\end{equation*}
	such that, as $k\to\iy$,
	\begin{equation}\label{11.23.22}
		\begin{split}
			\vp_{\de_k} &\ra \wt\vp_\be  \text{ weakly in } L^2(0,T; H^2(\Om_\be;A_\be)\cap H_{\Ga_\be}^1(\Om_\be;A_\be)),\\
			\pt_t\vp_{\de_k} &\ra \pt_t\wt\vp_\be  \text{ weakly in } L^2(Q_\be),\\
			\vp_{\de_k} &\ra \wt\vp_\be  \text{ strongly in } L^2(Q_\be) \text{ and } C([0,T]; L^2(\Om_\be)).
		\end{split}
	\end{equation}
	Since $\pt_t\vp_{\de_k} - \Div(A_\be\nabla\vp_{\de_k}) = f|_{Q_\be}$ a.e. in $Q_\be$ with $\vp_{\de_k}=0$ on $\Ga_\be\ts (0,T)$ and $\vp_{\de_k}(\cdot,0)=\vp^0|_{\Om_\be}$ in $L^2(\Om_\be)$ for all $\de_k\in (0,\be)$, passing to the limit as $k\to\iy$ via \eqref{11.23.22} yields
	\begin{equation}\label{11.23.23}
		\begin{cases}
			\pt_t\wt\vp_\be - \Div(A_\be\nabla\wt\vp_\be) = f|_{Q_\be} &\text{a.e. in } Q_\be,\\
			\wt\vp_\be = 0 &\text{on } \Ga_\be\ts (0,T),\\
			\wt\vp_\be(\cdot,0) = \vp^0|_{\Om_\be} &\text{in } L^2(\Om_\be).
		\end{cases}
	\end{equation}
	
By a diagonal argument across $\be_j=j^{-1}\de_0$ ($j\in\N^*$), we extract a common subsequence (still denoted by $\{\de_k\}$) and obtain a global function $\wt\vp$ on $Q$ such that $\wt\vp|_{Q_\be}=\wt\vp_\be$ for every $\be\in (0,\de_0]$. By weak lower semicontinuity together with \eqref{11.23.21} and \eqref{11.23.23}, we have
	\begin{equation*}
		\wt\vp \in L^2(0,T; H^2(\Om;A)\cap H_\Ga^1(\Om;A)) \cap H^1(0,T; L^2(\Om)),
	\end{equation*}
	and $\wt\vp$ satisfies $\pt_t\wt\vp - \Div(A\nabla\wt\vp) = f$ a.e. in $Q$ with $\wt\vp(\cdot,0)=\vp^0$ in $L^2(\Om)$. Since $\al\in [1,2)$, it follows from Theorem \ref{11.23.T1} that $H^2(\Om;A)\cap H_\Ga^1(\Om;A) = D(\mcA)$ (so that $\wt\vp$ automatically satisfies the homogeneous conormal boundary condition $\f{\pt\wt\vp}{\pt\nu_A}=0$ on the degenerate boundary). Consequently, $\wt\vp$ is a strong (and hence weak) solution to \eqref{11.23.2} at $\de=0$ associated with $(\vp^0, f)$. By the uniqueness of weak solutions established in Theorem \ref{11.23.T2}, we conclude that $\wt\vp=\vp_0$ in $Q$. Since the limit $\vp_0$ is independent of the choice of the subsequence $\{\de_k\}$, the entire family $\{\vp_\de\}$ converges to $\vp_0$ as $\de\to 0^+$ in the sense of \eqref{11.23.19}.

	\textit{Step 2}. Proof of \eqref{11.23.20}.
	
	Now assume additionally that $f=0$. By Theorem \ref{11.23.T2} (estimate \eqref{11.23.5} with $f_\de=0$)  and the inclusion $\Om_{4\de_0}\s \Om_\de$, we obtain
	\begin{equation}\label{11.23.24}
		\|\vp_\de\|_{L^2(0,T; H^2(\Om_{4\de_0}))} + \|\pt_t\vp_\de\|_{L^2(0,T; L^2(\Om_{4\de_0}))} \leq C\|\vp^0\|_{H_\Ga^1(\Om;A)},
	\end{equation}
	where $C>0$ is independent of $\de\in [0,\de_0]$.
	
	Fix any fractional exponent $s\in (\f{3}{2}, 2)$. By the compact embedding $H^2(\Om_{4\de_0})\hookrightarrow\hookrightarrow H^s(\Om_{4\de_0})\hookrightarrow L^2(\Om_{4\de_0})$, the Aubin-Lions lemma, and \eqref{11.23.19}, we deduce from \eqref{11.23.24} that
	\begin{equation}\label{11.23.25}
		\vp_\de \ra \vp_0 \text{ strongly in } L^2(0,T; H^s(\Om_{4\de_0}))  \text{ as } \de\to 0^+.
	\end{equation}
	Finally, since $s>\f{3}{2}$, the Sobolev trace theorem ensures that the normal derivative trace $\psi \mapsto \f{\pt\psi}{\pt\nu}|_{\Ga_{4\de_0}}$ is a bounded linear operator from $H^s(\Om_{4\de_0})$ into $H^{s-\f{3}{2}}(\Ga_{4\de_0})\hookrightarrow L^2(\Ga_{4\de_0})$. Combining this continuity with \eqref{11.23.25} gives
	\begin{equation*}
		\left\|\f{\pt\vp_\de}{\pt\nu} - \f{\pt\vp_0}{\pt\nu}\right\|_{L^2(0,T; L^2(\Ga_{4\de_0}))} \leq C\|\vp_\de - \vp_0\|_{L^2(0,T; H^s(\Om_{4\de_0}))} \ra 0  \text{ as } \de\to 0^+,
	\end{equation*}
	which proves \eqref{11.23.20} and completes the proof of Theorem \ref{11.23.T3}.
\end{proof}

\begin{remark}
	When we take $\Om_\de=\Om\cap \{x=(x_1,x_2)\in\R^2\colon x_2>\de\}$ with $\de\in (0,\de_0)$, all the results in this section still hold. Moreover, the observability of 
	\begin{equation}\label{11.25.15}
		\begin{cases}
			\pt_t\vp +\Div(A\nabla\vp)=0, &\mbox{in }Q, \\
			\vp=0, &\mbox{on }\Ga\ts (0,T), \\
			\f{\pt\vp}{\pt \nu_A}=0, &\mbox{on }\Ga_2^0\ts (0,T), \\
			\vp(T)=\vp^T, &\mbox{in }\Om
		\end{cases}
	\end{equation}
	can be reduced to the observability of the shape design equations
	\begin{equation*}
		\begin{cases}
			\pt_t\vp_\de +\Div(A_\de\nabla\vp_\de)=0, &\mbox{in }Q_\de, \\
			\vp_\de=0, &\mbox{on }\Ga_\de\ts (0,T), \\
			\f{\pt\vp_\de}{\pt \nu_{A_\de}}=0, &\mbox{on }\Ga_2^\de\ts (0,T), \\
			\vp_\de(T)=\vp_\de^T, &\mbox{in }\Om_\de.
		\end{cases}
	\end{equation*}
	However, deriving a Carleman estimate on $\Om_\de$ inevitably involves boundary terms at the corners $\pt\Ga_2^\de$ of $\Om_\de$, which are difficult to estimate. Therefore, we shall establish a different Carleman estimate directly for equation \eqref{11.25.15} (via the regularized system \eqref{11.25.3}).
\end{remark}

\section{Carleman estimate and observability}\label{S4}

In this section, we establish the observability inequality for equation \eqref{11.25.15}.

Let $\al\in [1,2)$ and $\rho\in (0,1)$. Set
\begin{equation}\label{11.25.1}
	\Theta =\f{1}{[t(T-t)]^4}, \  r=x_2+\rho,\  A^\rho=\diag(1,r^\al), \  \eta_\rho=r^{2-\al}-\rho^{2-\al}, \  \xi_\rho=\Theta (\ga-\eta_\rho),
\end{equation}
where $\ga=\max_{x\in \Om}x_2^{2-\al}+1$. Then the concavity of $x_2^{2-\al}$ gives $0\leq \eta_\rho\leq x_2^{2-\al}$, hence $1\leq \ga-\eta_\rho\leq \ga$. Define
\begin{equation}\label{11.25.2}
	\mcA^\rho u =-\Div(A^\rho \nabla u) \mbox{ on }\Om,\quad a^\rho(u,v)=\int_\Om \nabla u\cdot A^\rho \nabla v\df x.
\end{equation}
Denote
\begin{equation*}
	\Ga_+=\left\{x\in\pt\Om\colon \nu_2(x)>0\right\}.
\end{equation*}

\subsection{An approximation lemma}\label{S4.1}

\begin{lemma}\label{11.25.L1}
	Let $\vp^T\in H_\Ga^1(\Om;A)$ and $f\in L^2(Q)$. There exists a family $\{\vp_\rho^T\}_{\rho\in(0,1)}\subset H_0^1(\Om)$ satisfying
	\begin{equation*}
		\vp_\rho^T\ra \vp^T\mbox{ in } L^2(\Om),  \mbox{ and }  a^\rho(\vp_\rho^T,\vp_\rho^T)\ra \|\vp^T\|_{H_0^1(\Om;A)}^2  \mbox{ as }\rho\ra 0^+,
	\end{equation*}
	such that the weak solutions $\vp^\rho\ (0< \rho<1)$ of
	\begin{equation}\label{11.25.3}
		\begin{cases}
			\pt_t\vp^\rho +\Div(A^\rho\nabla\vp^\rho)=f, &\mbox{in }Q, \\
			\vp^\rho=0, &\mbox{on }\Sigma, \\
			\vp^\rho(T)=\vp_\rho^T, &\mbox{in }\Om,
		\end{cases}
	\end{equation}
	satisfy, as $\rho\ra 0^+$,
	\begin{equation}\label{11.25.4}
		\begin{split}
			\vp^\rho &\ra \vp \mbox{ strongly in }C([0,T]; L^2(\Om)), \\
			\f{\pt \vp^\rho}{\pt\nu} &\ra \f{\pt \vp}{\pt \nu} \mbox{ strongly in }L^2(\Ga_+\ts (0,T)),
		\end{split}
	\end{equation}
	where $\vp$ is the weak solution of
	\begin{equation}\label{11.25.5}
		\begin{cases}
			\pt_t\vp +\Div(A\nabla\vp)=f, &\mbox{in }Q, \\
			\vp=0, &\mbox{on }\Ga\ts (0,T), \\
			\f{\pt\vp}{\pt \nu_A}=0, &\mbox{on }\Ga_2^0\ts (0,T), \\
			\vp(T)=\vp^T, &\mbox{in }\Om.
		\end{cases}
	\end{equation}
\end{lemma}

\begin{proof}
	We prove this lemma in three steps. 
	
	{\it Step 1.} Construction of the terminal data $\vp_\rho^T$.
	
	By Lemma \ref{11.24.L1}, there exists a sequence $\{\vp_j^T\}_{j\in\N^*}\s C_0^\iy(\Om)$ such that
	\begin{equation*}
		\vp_j^T \ra \vp^T  \mbox{ strongly in } H_0^1(\Om;A) \mbox{ (and hence in } L^2(\Om)\mbox{)}  \mbox{ as } j\ra\iy.
	\end{equation*}
	For each fixed $j\in \N^*$, using the elementary inequality $(x_2+\rho)^\al - x_2^\al \le \al \rho (x_2+\rho)^{\al-1}$ for $\al\in [1,2)$, we obtain
	\begin{equation*}
		\begin{split} 
		0 
		&\le a^\rho(\vp_j^T,\vp_j^T) - \|\vp_j^T\|_{H_0^1(\Om;A)}^2\\ 
		&= \int_\Om \big((x_2+\rho)^\al - x_2^\al\big) |\pt_{x_2}\vp_j^T|^2 \df x 
		\le \al \rho \sup_{x\in\Om}(x_2+1)^{\al-1} \|\nabla \vp_j^T\|_{L^2(\Om)}^2 \ra 0
		\end{split} 
	\end{equation*}
	as $\rho\ra 0^+$. Consequently, for every $j\in \N^*$, we can choose $\rho_j \in (0, j^{-1})$ strictly decreasing with respect to $j$ such that for all $\rho \in (0, \rho_j]$,
	\begin{equation*}
		\big| a^\rho(\vp_j^T,\vp_j^T) - \|\vp_j^T\|_{H_0^1(\Om;A)}^2 \big| < \f{1}{j}.
	\end{equation*}
If we set $\vp_\rho^T := \vp_j^T$ for $\rho \in (\rho_{j+1}, \rho_j]$ ($j\in \N^*$, and $\vp_\rho^T := \vp_1^T$ for $\rho \in (\rho_1, 1)$), it follows immediately from the triangle inequality that $\vp_\rho^T \in C_0^\iy(\Om)\s H_0^1(\Om)$ satisfies
	\begin{equation}\label{11.25.9}
		\vp_\rho^T \ra \vp^T \mbox{ in } L^2(\Om), \mbox{ and }  a^\rho(\vp_\rho^T,\vp_\rho^T) \ra \|\vp^T\|_{H_0^1(\Om;A)}^2  \mbox{ as } \rho\ra 0^+.
	\end{equation}

	{\it Step 2.} Energy estimates and interior convergence $\vp^\rho \ra \vp$.
	
	For each $\rho\in (0,1)$, since $H_0^1(\Om)\hra H_\Ga^1(\Om;A)$ is continuous by Lemma \ref{11.24.L1}, it follows from Remark \ref{11.22.R1} that
	\begin{equation}\label{11.25.7}
		\int_\Om u^2\df x\leq C_P\int_\Om \nabla u\cdot A\nabla u\df x\leq C_P a^\rho(u,u) \quad \mbox{for all } u\in H_0^1(\Om).
	\end{equation}
	
	For each $\rho\in (0,1)$, since \eqref{11.25.3} is uniformly parabolic, multiplying \eqref{11.25.3} by $\vp^\rho$, integrating over $\Om\ts (t,T)$, and applying \eqref{11.25.7} together with Young's inequality, we obtain
	\begin{equation}\label{11.25.6}
		\|\vp^\rho\|_{L^\iy(0,T; L^2(\Om))}^2+\int_0^T a^\rho(\vp^\rho,\vp^\rho)\df t\leq C\left(\|\vp_\rho^T\|_{L^2(\Om)}^2+\|f\|_{L^2(Q)}^2\right),
	\end{equation}
	where the positive constant $C$ depends only on $\Om$ and $\al$. Next, multiplying \eqref{11.25.3} by $\pt_t\vp^\rho$ and integrating over $\Om\ts (t,T)$, and also using the relation $\mcA^\rho\vp^\rho=\pt_t\vp^\rho-f$, we deduce that
	\begin{equation}\label{11.25.8}
		\esssup_{t\in [0,T]} a^\rho(\vp^\rho(t),\vp^\rho(t))+\|\pt_t\vp^\rho\|_{L^2(Q)}^2+\|\mcA^\rho \vp^\rho\|_{L^2(Q)}^2\leq 3\left(a^\rho(\vp_\rho^T,\vp_\rho^T)+\|f\|_{L^2(Q)}^2\right).
	\end{equation}
	
By \eqref{11.25.9}, \eqref{11.25.6},  \eqref{11.25.8}, and noting that $A^\rho\geq A$ in $\Om$, there exist a subsequence of $\{\vp^\rho\}_{\rho\in (0,1)}$ (still denoted by $\{\vp^\rho\}$) and a function
	\begin{equation}\label{11.25.11}
		\wt\vp\in L^\iy(0,T; H_\Ga^1(\Om;A)) \mbox{ with }  \pt_t\wt\vp\in L^2(Q), 
	\end{equation}
	such that 
	\begin{equation}\label{11.25.10}
		\begin{split}
			\vp^\rho &\ra \wt\vp  \mbox{ weakly-* in }L^\iy(0,T; H_\Ga^1(\Om;A)), \\
			\pt_t\vp^\rho &\ra \pt_t\wt\vp  \mbox{ weakly in } L^2(Q).
		\end{split}
	\end{equation}
	Moreover, since the embedding $H_{\Ga}^1(\Om;A)\hra L^2(\Om)$ is compact, the Aubin-Lions lemma yields
	\begin{equation}\label{11.25.12}
		\vp^\rho \ra \wt\vp  \mbox{ strongly in } C([0,T]; L^2(\Om)). 
	\end{equation}
	
	For each $\psi\in C^\iy(\ol Q)$ satisfying $\dist(\supp \psi(t), \pt\Om)>0$ for all $t\in [0,T]$ and $\psi(0)=0$, multiplying \eqref{11.25.3} by $\psi$ and integrating over $Q$, we obtain 
	\begin{equation*}
		\int_\Om \vp_\rho^T\psi(T)\df x-\iint_Q \vp^\rho \pt_t\psi\df x\df t-\iint_Q \nabla\vp^\rho\cdot A^\rho \nabla\psi\df x\df t=\iint_Q f\psi\df x\df t. 
	\end{equation*}
	Passing to the limit along the subsequence and using \eqref{11.25.9}, \eqref{11.25.10}, \eqref{11.25.12}, and the fact that $A^\rho\ra A$ uniformly in $\Om$, we get
	\begin{equation}\label{11.25.13}
		\int_\Om \vp^T\psi(T)\df x-\iint_Q \wt \vp \pt_t\psi\df x\df t-\iint_Q \nabla\wt \vp\cdot A \nabla\psi\df x\df t=\iint_Q f\psi\df x\df t. 
	\end{equation}
	This implies that 
	\begin{equation*}
		\pt_t\wt\vp+\Div(A\nabla\wt\vp)=f \mbox{ in }Q \mbox{ in the sense of distributions},  \mbox{ and }  \wt\vp(T)=\vp^T. 
	\end{equation*}
	Note that by Lemma \ref{11.24.L1} (i.e., the density of $C_0^\iy(\Om)$ in $H_\Ga^1(\Om;A)$), \eqref{11.25.13} also holds for all $\psi\in C^\iy(\ol Q)$ with $\dist(\supp \psi(t), \Ga)>0$ for all $t\in [0,T]$ and $\psi(0)=0$. Combining this with Proposition \ref{11.23.P1} (for $\de=0$), we deduce that $\wt\vp$ is a weak solution of \eqref{11.25.5}. By Theorem \ref{11.23.T2} (uniqueness of the weak solution), it follows that $\wt\vp=\vp$. Consequently, the whole family $\vp^\rho$ converges to $\vp$ strongly in $C([0,T]; L^2(\Om))$, and hence in $L^2(Q)$, as $\rho\ra 0^+$.
	
	{\it Step 3.} Boundary convergence of the normal derivative.

	By \eqref{11.22.3}, $\nu_2<0$ on the boundary with
	$0\leq x_2<8\delta_0$. Thus
	$\overline{\Gamma_+}\subset\partial\Omega\cap\{x_2\geq8\delta_0\}$.
	Use the fixed domain $D=\Omega_{\delta_0}$ and the cutoff
	$\zeta$ from the proof of Theorem~\ref{11.23.T2}.
	Set $w^\rho=\zeta\varphi^\rho|_D$.
	Then $w^\rho\in L^2(0,T;H_0^1(D))$ and
	\[
	\begin{split}
		-\operatorname{div}(A^\rho\nabla w^\rho)
		&=\zeta(\partial_t\varphi^\rho-f)
		-2(x_2+\rho)^\alpha\zeta'\partial_{x_2}\varphi^\rho\\
		&\quad-\big((x_2+\rho)^\alpha\zeta'\big)'\varphi^\rho.
	\end{split}
	\]
	The matrices $A^\rho|_D$ have ellipticity and $C^1$ bounds
	independent of $\rho\in(0,1)$. The boundary $H^2$ estimate,
	\eqref{11.25.8}, and \eqref{11.25.9} give
	\[
	\sup_{0<\rho<\rho_1}
	\|w^\rho\|_{L^2(0,T;H^2(D))}<\infty
	\]
	for some $\rho_1>0$. By \eqref{11.25.12},
	$w^\rho\to w:=\zeta\varphi|_D$ in $L^2(D\times(0,T))$,
	and $w\in L^2(0,T;H^2(D))$ by weak compactness.
	For every $\sigma\in(3/2,2)$, interpolation yields
	\[
	\|w^\rho-w\|_{L^2(0,T;H^\sigma(D))}
	\leq C\|w^\rho-w\|_{L^2(0,T;H^2(D))}^{\sigma/2}
	\|w^\rho-w\|_{L^2(D\times(0,T))}^{1-\sigma/2}
	\longrightarrow0.
	\]
	Since $\zeta=1$ in a neighborhood of $\overline{\Gamma_+}$,
	the normal derivative trace theorem gives
	\[
	\partial_\nu\varphi^\rho\longrightarrow\partial_\nu\varphi
	\quad\text{in }L^2(\Gamma_+\times(0,T)).
	\]
	This completes the proof.
\end{proof}

\subsection{Carleman estimate}\label{S4.2}

Let $\vp^T\in H_\Ga^1(\Om;A)$ and $f\in L^2(Q)$. Let $\rho\in (0,1)$ and $s>0$. Define
\begin{equation*}
	z=\vp^\rho e^{-s\xi_\rho} \iff \vp^\rho=e^{s\xi_\rho}z. 
\end{equation*}
Then, 
\begin{equation}\label{11.25.14}
	\begin{split}
		&(i)\ z(t)=0 \mbox{ and } \nabla z(t)=0 \mbox{ at }t=0,T, \\
		&(ii)\ z=\pt_t z=0 \mbox{ on }\Sigma, 
	\end{split}
\end{equation}
and \eqref{11.25.3} can be rewritten as
\begin{equation*}
	e^{-s\xi_\rho}f=P_1z+P_2z, 
\end{equation*}
where 
\begin{equation*}
	\begin{split}
		P_1z &= \pt_t z + 2s\nabla z \cdot A^\rho\nabla\xi_\rho + sz\Div(A^\rho \nabla\xi_\rho), \\
		P_2z &= \Div(A^\rho\nabla z) + sz\pt_t\xi_\rho + s^2 z \nabla\xi_\rho \cdot A^\rho\nabla\xi_\rho. 
	\end{split}
\end{equation*}

Since \eqref{11.25.3} is a uniformly parabolic equation with homogeneous Dirichlet boundary conditions, all the integrations by parts performed below are justified.

\vs{3mm}

\noindent\textbf{1. Computation}

\vs{3mm}

(1) Writing $A^\rho=\diag(A_1^\rho,A_2^\rho)$ with $A_1^\rho=1$ and $A_2^\rho=r^\al$, from \eqref{11.25.14} and integration by parts, we obtain 
\begin{equation*}
	\begin{split}
		&(P_1z,P_2z)_{L^2(Q)}\\
		&=2s\iint_\Sigma  (\nabla z\cdot A^\rho \nu) (\nabla z\cdot A^\rho\nabla\xi_\rho)\df S\df t-s\iint_\Sigma  (\nabla z\cdot A^\rho\nabla z)(\nabla\xi_\rho\cdot A^\rho\nu)\df S\df t\\
		&\hspace{4.5mm}-2s\iint_Q A^\rho \nabla z\cdot \left[(D^2\xi_\rho)A^\rho \nabla z\right]\df x\df t-s\iint_Q zA^\rho\nabla z\cdot \nabla [\Div(A^\rho\nabla\xi_\rho)]\df x\df t\\
		&\hspace{4.5mm}-2s\sum_{i,j=1}^2\iint_Q A_i^\rho\f{\pt z}{\pt x_i} \f{\pt A_j^\rho}{\pt x_i}\f{\pt z}{\pt x_j}\f{\pt \xi_\rho}{\pt x_j}\df x\df t+s\sum_{i,j=1}^2\iint_Q \left(\f{\pt z}{\pt x_i}\right)^2\f{\pt A_i^\rho }{\pt x_j}A_j^\rho\f{\pt \xi_\rho}{\pt x_j}\df x\df t\\
		&\hspace{4.5mm}-\f{s}{2}\iint_Q z^2\pt_{tt}\xi_\rho\df x\df t-2s^2\iint_Q z^2\nabla\xi_\rho\cdot A^\rho\nabla  \pt_t\xi_\rho\df x\df t\\
		&\hspace{4.5mm}-s^3\iint_Q z^2A^\rho\nabla\xi_\rho\cdot \nabla(\nabla\xi_\rho\cdot A^\rho\nabla\xi_\rho)\df x\df t. 
	\end{split}
\end{equation*}

\noindent\textbf{2. Estimates}

\vs{3mm}

(2) Recall that $\Theta(t)=[t(T-t)]^{-4}$ for $t\in(0,T)$. Direct computations yield
\begin{equation*} 
	|\Theta'|\leq 4T\Theta^\f{5}{4},\quad |\Theta\Theta'|\leq 4T\Theta^\f{9}{4},\quad |\Theta''|\leq 20T^2\Theta^\f{3}{2}, 
\end{equation*}
and, denoting $e_1=(1,0)$ and $e_2=(0,1)$,
\begin{equation*}
	\nabla \eta_\rho=(2-\al)r^{1-\al}e_2, \quad \nabla\xi_\rho=-(2-\al)\Theta r^{1-\al}e_2,\quad \nabla\pt_t\xi_\rho=-(2-\al)\Theta' r^{1-\al}e_2, 
\end{equation*}
and 
\begin{equation*}
	\f{\pt A_i^\rho}{\pt x_j}=\al r^{\al-1}\de_{i2}\de_{j2},\quad A^\rho \nabla\xi_\rho=-(2-\al)\Theta re_2, \quad \Div(A^\rho\nabla\xi_\rho)=-(2-\al)\Theta, 
\end{equation*}
and 
\begin{equation*}
	\f{\pt^2 \xi_\rho}{\pt x_i\pt x_j}=-(2-\al)(1-\al)\Theta r^{-\al}\de_{i2}\de_{j2}, \quad \pt_{tt}\xi_\rho=\Theta''(\ga-\eta_\rho),\quad \nabla [\Div(A^\rho\nabla\xi_\rho)]=0, 
\end{equation*}
and 
\begin{equation*}
	\nabla\xi_\rho\cdot A^\rho\nabla\pt_t\xi_\rho=(2-\al)^2\Theta\Theta'r^{2-\al},\quad A^\rho\nabla\xi_\rho\cdot \nabla (\nabla\xi_\rho\cdot A^\rho\nabla\xi_\rho)=-(2-\al)^4\Theta^3r^{2-\al}, 
\end{equation*}
where $\de_{ij}=1$ for $i=j$ and $\de_{ij}=0$ for $i\neq j$. 

From \eqref{11.25.14}, we have $\nabla z=\f{\pt z}{\pt \nu}\nu$ on $\Sigma$, and hence
\begin{equation*}
	(\nabla z\cdot A^\rho\nu)(\nabla z\cdot A^\rho\nabla\xi_\rho)=(\nabla z\cdot A^\rho\nabla z)(\nabla\xi_\rho\cdot A^\rho\nu)=-(2-\al)\Theta r\left(\f{\pt z}{\pt \nu}\right)^2 (\nu\cdot A^\rho\nu)\nu_2. 
\end{equation*}

(3) Combining (1) and (2), we obtain 
\begin{equation*}
	\begin{split}
		(P_1z,P_2z)_{L^2(Q)}
		&\geq -(2-\al)s\iint_\Sigma \Theta \left(\f{\pt z}{\pt \nu}\right)^2(\nu \cdot A^\rho\nu)(r\nu_2)\df S\df t\\
		&\hspace{4.5mm}+(2-\al)^2s\iint_Q \Theta r^\al (\pt_{x_2}z)^2\df x\df t+(2-\al)^4s^3\iint_Q \Theta^3z^2r^{2-\al}\df x\df t\\
		&\hspace{4.5mm}-\f{s}{2}\iint_Q\Theta''  z^2(\ga-\eta_\rho)\df x\df t-2(2-\al)^2s^2\iint_Q \Theta\Theta' z^2r^{2-\al}\df x\df t. 
	\end{split}
\end{equation*}

By a direct computation, we obtain 
\begin{equation*}
	s^2\iint_Q |\Theta \Theta'|z^2r^{2-\al}\df x\df t\leq 4Ts^2\iint_Q \Theta^\f{9}{4} z^2r^{2-\al}\df x\df t\leq \f{1}{16}T^7s^2\iint_Q \Theta^3z^2r^{2-\al}\df x\df t. 
\end{equation*}

From \eqref{11.25.14}, integration by parts gives
\begin{equation}\label{11.25.16}
	0=\int_\Om \pt_{x_2}(r z^2)\df x=\int_\Om z^2\df x+2\int_\Om rz\pt_{x_2}z\df x.
\end{equation}
Using this identity together with $|\Theta''|\leq 20T^2\Theta^\f{3}{2}$, $0\leq \ga-\eta_\rho\leq \ga$, Young's inequality, and $\Theta(t)\geq 256T^{-8}$ for $t\in(0,T)$, we obtain
\begin{equation*}
	\begin{split}
		s\iint_Q |\Theta''|z^2(\ga-\eta_\rho)\df x\df t
		&\leq 20\ga T^2s\iint_Q \Theta^\f{3}{2} z^2\df x\df t=-40\ga T^2s\iint_Q \Theta^\f{3}{2} r z\pt_{x_2}z\df x\df t\\
		&\leq \ga^\f{2}{3}T^4\iint_Q \Theta r^\al |\pt_{x_2}z|^2\df x\df t+400\ga^\f{4}{3}s^2\iint_Q \Theta^2z^2r^{2-\al}\df x\df t\\
		&\leq \ga^\f{2}{3}T^4\iint_Q \Theta r^\al |\pt_{x_2}z|^2\df x\df t+\f{25}{16}\ga^\f{4}{3}T^8s^2\iint_Q \Theta^3z^2r^{2-\al}\df x\df t.
	\end{split}
\end{equation*}

Hence, there exists $\la_0\geq 1$, depending only on $\al$ and $\ga$, such that for all $\la\geq \la_0$ and all $s=\la (1+T^8)$, 
\begin{equation*}
	\begin{split}
		(P_1z,P_2z)_{L^2(Q)}
		&\geq -(2-\al)s\iint_\Sigma \Theta \left(\f{\pt z}{\pt \nu}\right)^2(\nu \cdot A^\rho\nu)(r\nu_2)\df S\df t\\
		&\hspace{4.5mm}+Cs\iint_Q \Theta r^\al (\pt_{x_2}z)^2\df x\df t+Cs^3\iint_Q \Theta^3z^2r^{2-\al}\df x\df t, 
	\end{split}
\end{equation*}
where the positive constant $C$ depends only on $\al$.

Finally, using the identity $0=\int_\Om \pt_{x_2}(rz^2)\df x$ and Young's inequality, we have
\begin{equation*}
	\begin{split}
		s^2\iint_Q \Theta^2 z^2\df x\df t
		&=-2s^2\iint_Q \Theta^2rz\pt_{x_2}z\df x\df t\\
		&\leq s\iint_Q \Theta r^\al|\pt_{x_2}z|^2\df x\df t+s^3\iint_Q \Theta^3 z^2r^{2-\al}\df x\df t. 
	\end{split}
\end{equation*}
Combining this with $\|e^{-s\xi_\rho}f\|_{L^2(Q)}^2=\|P_1z+P_2z\|_{L^2(Q)}^2\geq 2(P_1z,P_2z)_{L^2(Q)}$ and noting that $(\nu\cdot A^\rho\nu)(r\nu_2)\leq 0$ on $\Sigma-(\Ga_+\ts(0,T))$, we obtain, for all $s\geq s_0$ with $s_0=\la_0(1+T^8)$, 
\begin{equation}\label{11.25.17}
	\begin{split}
		&s\iint_Q \Theta r^\al (\pt_{x_2}z)^2\df x\df t+s^2\iint_Q \Theta^2z^2\df x\df t+s^3\iint_Q \Theta^3z^2r^{2-\al}\df x\df t\\
		&\leq C\|e^{-s\xi_\rho}f\|_{L^2(Q)}^2+Cs\iint_{\Ga_+\ts(0,T)}\Theta \left(\f{\pt z}{\pt \nu}\right)^2 (\nu\cdot A^\rho\nu)(r\nu_2)\df S\df t\\
		&\leq C\|e^{-s\xi_\rho}f\|_{L^2(Q)}^2+Cs\iint_{\Ga_+\ts (0,T)}\Theta \left(\f{\pt z}{\pt \nu}\right)^2 \df S\df t, 
	\end{split}
\end{equation}
where the positive constant $C$ depends only on $\al$ and $M$. 

(4) Since $\vp^\rho=0$ on $\Sigma$, we have $\f{\pt z}{\pt\nu}=e^{-s\xi_\rho}\f{\pt\vp^\rho}{\pt\nu}$ on $\Ga_+\ts(0,T)$. Moreover, noting that
\begin{equation*}
	\pt_{x_2}\vp^\rho=\pt_{x_2}\left(e^{s\xi_\rho}z\right)=e^{s\xi_\rho}\left(\pt_{x_2}z-(2-\al)s\Theta zr^{1-\al}\right), 
\end{equation*}
from \eqref{11.25.17} (and the triangle inequality), we obtain for all $s\geq s_0$, 
\begin{equation*}
	\begin{split}
		&s\iint_Q \Theta r^\al (\pt_{x_2}\vp^\rho)^2e^{-2s\xi_\rho}\df x\df t+s^2\iint_Q \Theta^2(\vp^\rho)^2e^{-2s\xi_\rho}\df x\df t+s^3\iint_Q \Theta^3(\vp^\rho)^2r^{2-\al}e^{-2s\xi_\rho}\df x\df t\\
		&\leq C\|e^{-s\xi_\rho}f\|_{L^2(Q)}^2+Cs\iint_{\Ga_+\ts (0,T)}\Theta \left(\f{\pt \vp^\rho}{\pt \nu}\right)^2 e^{-2s\xi_\rho}\df S\df t, 
	\end{split}
\end{equation*}
where the positive constant $C$ depends only on $\al$ and $M$. 

Now let $f=0$. By the energy identity for \eqref{11.25.3}, we have 
\begin{equation*}
	\|\vp^\rho(t_1)\|_{L^2(\Om)}\leq \|\vp^\rho(t_2)\|_{L^2(\Om)} \quad \mbox{for all } 0\leq t_1\leq t_2\leq T. 
\end{equation*}
Note that for all $(x,t)\in \Om\ts \left[\f{T}{4},\f{3T}{4}\right]$, we have $t(T-t)\geq \f{3}{16}T^2$, and hence
\begin{equation*}
	\Theta(t)\leq \left(\f{16}{3T^2}\right)^4\leq 810 T^{-8}, \quad 0\leq \xi_\rho(x,t)\leq \ga\Theta(t)\leq 810\ga T^{-8}.
\end{equation*}
Taking $s=s_0=\la_0(1+T^8)$, we deduce that
\begin{equation*}
	\begin{split}
		s_0^2\int_\Om \left|\vp^\rho(0)\right|^2\df x
		&\leq \f{2}{T}s_0^2\int_\f{T}{4}^\f{3T}{4}\int_\Om |\vp^\rho(t)|^2\df x\df t\\
		&\leq C T^{15}e^{C s_0 T^{-8}}s_0^2\iint_Q \Theta^2 (\vp^\rho)^2e^{-2s_0\xi_\rho}\df x\df t\\
		&\leq C T^{15}e^{C s_0 T^{-8}}s_0\iint_{\Ga_+\ts (0,T)}\Theta \left(\f{\pt \vp^\rho}{\pt \nu}\right)^2 e^{-2s_0\xi_\rho}\df S\df t\\
		&\leq C T^{15}e^{C(1+T^{-8})}s_0\iint_{\Ga_+\ts (0,T)} \left(\f{\pt \vp^\rho}{\pt \nu}\right)^2 \df S\df t,
	\end{split}
\end{equation*}
where in the last step we used the uniform boundedness of $\Theta e^{-2s_0\xi_\rho}$ on $\Ga_+\ts(0,T)$ (since $\ga-\eta_\rho\geq 1$, so $\Theta e^{-2s_0\xi_\rho}\leq \Theta e^{-2s_0\Theta}\leq (2s_0 e)^{-1}$).

\subsection{Observability inequality}\label{S4.3}

(5) Passing to the limit as $\rho\ra 0^+$ in (4) and using \eqref{11.25.4} together with \eqref{11.25.12}, we obtain 
\begin{equation}\label{11.26.1}
	\int_\Om |\vp(0)|^2\df x\leq C\iint_{\Ga_+\ts (0,T)}\left(\f{\pt \vp}{\pt \nu}\right)^2\df S\df t, 
\end{equation}
where the positive constant $C$ depends only on $\al$, $M$, and $T$. 

(6) 
Let $\varphi^T\in L^2(\Omega)$ and $f=0$,
and put $\tau=T/2$. Write
$c_n=(\varphi^T,\Phi_n^0)_{L^2(\Omega)}$. By the spectral
representation of the solution,
\[
\varphi(t)=\sum_{n=1}^{\infty}
c_n e^{-\lambda_n^0(T-t)}\Phi_n^0.
\]
Consequently,
\[
\begin{split}
	\|\varphi(\tau)\|_{H_\Gamma^1(\Omega;A)}^2
	&=\sum_{n=1}^{\infty}
	\lambda_n^0e^{-\lambda_n^0T}|c_n|^2 \leq\frac{1}{eT}\|\varphi^T\|_{L^2(\Omega)}^2.
\end{split}
\]
Thus $\varphi|_{(0,\tau)}$ has terminal datum
$\varphi(\tau)\in H_\Gamma^1(\Omega;A)$.
Applying \eqref{11.26.1} with $T$ replaced by $\tau$ gives
\[
\|\varphi(0)\|_{L^2(\Omega)}^2
\leq C_{\alpha,M,T}\int_0^{T/2}\int_{\Gamma_+}
|\partial_\nu\varphi|^2\,\mathrm dS\,\mathrm dt.
\]
The integral is finite by the graph-norm energy estimate
for $H_\Gamma^1(\Omega;A)$ terminal data, followed by the
fixed-domain boundary localization used in the proof of
Lemma~\ref{11.25.L1}. This proves the following theorem.

\begin{theorem}\label{11.26.T1}
	Let $\alpha\in[1,2)$, $\varphi^T\in L^2(\Omega)$, and $f=0$.
	The weak solution of \eqref{11.25.15} satisfies
	$\partial_\nu\varphi\in L^2(\Gamma_+\times(0,T/2))$ and
	\[
	\|\varphi(0)\|_{L^2(\Omega)}^2
	\leq C_{\alpha,M,T}\int_0^{T/2}\int_{\Gamma_+}
	|\partial_\nu\varphi|^2\,\mathrm dS\,\mathrm dt.
	\]
	In particular, the same inequality holds with observation
	on $(0,T)$, where for general $L^2$ terminal data the
	full-time integral is understood as a nonnegative Lebesgue
	integral, possibly equal to $+\infty$.
\end{theorem}

\section{Appendix}\label{S5}

In this section, we provide the proofs of Theorem \ref{11.23.T2} and Proposition \ref{11.23.P1}.

\begin{proof}[Proof of Theorem \ref{11.23.T2}]
	We prove this theorem by the following steps. 
	
	\textit{Step 1}. Galerkin approximation.
	
Let $\{\Phi_n^\de\}_{n=1}^\iy$ be the eigenfunctions of $\mcA_\de$ associated with the eigenvalues $\{\la_n^\de\}_{n=1}^\iy$. These eigenfunctions form an orthonormal basis of $L^2(\Om_\de)$ and an orthogonal basis of $H_{\Ga_\de}^1(\Om_\de;A_\de)$. For each $m\in\N^*$, we define the approximate solution
	\begin{equation}\label{11.23.6}
		\vp_\de^m(x,t)=\sum_{n=1}^m \vp_{\de,n}^m(t)\Phi_n^\de(x),\quad f_\de^m(x,t)=\sum_{n=1}^m f_{\de,n}^m(t)\Phi_n^\de(x), 
	\end{equation}
	where the coefficients $\vp_{\de,n}^m(t)$ ($n=1,\cdots, m$) solve the system of ordinary differential equations
	\begin{equation}\label{11.23.7}
		\f{\df}{\df t}\vp_{\de,n}^m(t)+\la_n^\de \vp_{\de,n}^m(t)=f_{\de,n}^m(t), \  n=1,\cdots, m, 
	\end{equation}
	subject to the initial conditions
	\begin{equation}\label{11.23.8}
		\vp_{\de,n}^m(0)=(\vp_\de^0,\Phi_n^\de)_{L^2(\Om_\de)}, \quad f_{\de,n}^m=(f_\de,\Phi_n^\de)_{L^2(\Om_\de)}, \quad   n=1,\cdots, m.  
	\end{equation}
	Clearly, the initial value problem \eqref{11.23.7}-\eqref{11.23.8} admits a unique solution $\vp_{\de,n}^m\in H^1(0,T)$ ($n=1,\cdots, m$). Moreover, using the orthogonality of $\{\Phi_n^\de\}_{n=1}^\iy$, equation \eqref{11.23.7} is equivalent to
	\begin{equation}\label{11.23.9}
		(\pt_t\vp_\de^m, \Phi_n^\de)_{L^2(\Om_\de)}+(\vp_\de^m, \Phi_n^\de)_{H_{\Ga_\de}^1(\Om_\de;A_\de)}=(f_\de^m, \Phi_n^\de)_{L^2(\Om_\de)}, \ n=1,\cdots, m. 
	\end{equation}
	
	\textit{Step 2}. Energy estimates. 
	
	Multiplying both sides of \eqref{11.23.9} by $\vp_{\de,n}^m(t)$ and summing over $n=1,\cdots, m$, we obtain 
	\begin{equation*}
		(\pt_t\vp_\de^m, \vp_\de^m)_{L^2(\Om_\de)}+(\vp_\de^m, \vp_\de^m)_{H_{\Ga_\de}^1(\Om_\de;A_\de)}=(f_\de^m,\vp_\de^m)_{L^2(\Om_\de)}. 
	\end{equation*}
	Since $(\pt_t\vp_\de^m, \vp_\de^m)_{L^2(\Om_\de)}=\f{1}{2}\f{\df}{\df t}\|\vp_\de^m\|_{L^2(\Om_\de)}^2$, applying Young's inequality gives 
	\begin{equation}\label{11.23.10}
		\begin{split}
			\f{1}{2}\f{\df}{\df t}\|\vp_\de^m\|_{L^2(\Om_\de)}^2+\|\vp_\de^m\|_{H_{\Ga_\de}^1(\Om_\de;A_\de)}^2\leq \f{1}{2}\|f_\de^m\|_{L^2(\Om_\de)}^2+\f{1}{2}\|\vp_\de^m\|_{L^2(\Om_\de)}^2, 
		\end{split}
	\end{equation}
	and hence 
	\begin{equation*}
		\begin{split}
			\|\vp_\de^m(t)\|_{L^2(\Om_\de)}^2\leq e^t\left(\|\vp_\de^m(0)\|_{L^2(\Om_\de)}^2+\int_0^t\|f_\de^m(s)\|_{L^2(\Om_\de)}^2\,\df s\right), \text{ for all } t\in [0,T]
		\end{split}
	\end{equation*}
	by Gronwall's inequality \cite[Appendix B.2 (j), p.~708]{Evans}. Note that $\|\vp_\de^m(0)\|_{L^2(\Om_\de)}\leq \|\vp_\de^0\|_{L^2(\Om_\de)}$ and $\|f_\de^m(s)\|_{L^2(\Om_\de)}\leq \|f_\de(s)\|_{L^2(\Om_\de)}$. Combining this with \eqref{11.23.10} integrated over $(0,T)$, we deduce that 
	\begin{equation}\label{11.23.11}
		\begin{split}
			\|\vp_\de^m\|_{L^\iy(0,T; L^2(\Om_\de))}^2+2\|\vp_\de^m\|_{L^2(0,T; H_{\Ga_\de}^1(\Om_\de;A_\de))}^2\leq 2e^T\left(\|\vp_\de^0\|_{L^2(\Om_\de)}^2+\|f_\de\|_{L^2(Q_\de)}^2\right). 
		\end{split}
	\end{equation}
	
	Next, let $v\in H_{\Ga_\de}^1(\Om_\de;A_\de)$ with $\|v\|_{H_{\Ga_\de}^1(\Om_\de;A_\de)}\leq 1$. Write $v=v^1+v^2$ with $v^1\in \Span\{\Phi_n^\de\}_{n=1}^m$ and $(v^2,\Phi_n^\de)_{L^2(\Om_\de)}=0$ for all $n=1,\cdots, m$. Due to the orthogonality of $\{\Phi_n^\de\}_{n=1}^\iy$ in both $L^2(\Om_\de)$ and $H_{\Ga_\de}^1(\Om_\de;A_\de)$, we have $\|v^1\|_{H_{\Ga_\de}^1(\Om_\de;A_\de)}\leq \|v\|_{H_{\Ga_\de}^1(\Om_\de;A_\de)}\leq 1$. Multiplying both sides of \eqref{11.23.9} by $(v,\Phi_n^\de)_{L^2(\Om_\de)}$ and summing over $n=1,\cdots, m$, we get for a.e. $t\in (0,T)$,
	\begin{equation*}
		\begin{split}
			\lg \pt_t\vp_\de^m, v\rg_{H_{\Ga_\de}^{-1}(\Om_\de;A_\de), H_{\Ga_\de}^1(\Om_\de;A_\de)}
			&=\lg\pt_t\vp_\de^m, v^1\rg_{H_{\Ga_\de}^{-1}(\Om_\de;A_\de), H_{\Ga_\de}^1(\Om_\de;A_\de)}=(\pt_t\vp_\de^m, v^1)_{L^2(\Om_\de)}\\
			&=(f_\de,v^1)_{L^2(\Om_\de)}-(\vp_\de^m, v^1)_{H_{\Ga_\de}^1(\Om_\de;A_\de)}\\
			&\leq \|f_\de\|_{L^2(\Om_\de)}\|v^1\|_{L^2(\Om_\de)}+\|\vp_\de^m\|_{H_{\Ga_\de}^1(\Om_\de;A_\de)}\|v^1\|_{H_{\Ga_\de}^1(\Om_\de;A_\de)}\\
			&\leq \sqrt{C_{P}}\|f_\de\|_{L^2(\Om_\de)}+\|\vp_\de^m\|_{H_{\Ga_\de}^1(\Om_\de;A_\de)},
		\end{split}
	\end{equation*}
	where $C_P>0$ is the Poincar\'e constant independent of $\de$. Taking the supremum over all $v\in H_{\Ga_\de}^1(\Om_\de;A_\de)$ with $\|v\|_{H_{\Ga_\de}^1(\Om_\de;A_\de)}\leq 1$, squaring, integrating over $(0,T)$, and using \eqref{11.23.11}, we obtain 
	\begin{equation}\label{11.23.12}
		\begin{split}
			\|\pt_t\vp_\de^m\|_{L^2(0,T; H_{\Ga_\de}^{-1}(\Om_\de;A_\de))}\leq C\left(\|\vp_\de^0\|_{L^2(\Om_\de)}+\|f_\de\|_{L^2(Q_\de)}\right), 
		\end{split}
	\end{equation}
	where the positive constant $C$ depends only on $\al$, $M$, and $T$.
	
	\textit{Step 3}. Approximation, or convergence
	
	By the uniform estimates \eqref{11.23.11} and \eqref{11.23.12}, there exist a subsequence of $\{\vp_\de^m\}_{m=1}^\iy$ (still denoted by $\{\vp_\de^m\}_{m=1}^\iy$) and a function 
	\begin{equation*}
		\wt\vp_\de\in L^\iy(0,T; L^2(\Om_\de))\cap L^2(0,T; H_{\Ga_\de}^1(\Om_\de;A_\de)), \text{ with }  \pt_t\wt \vp_\de\in L^2(0,T; H_{\Ga_\de}^{-1}(\Om_\de;A_\de))
	\end{equation*}
	such that, as $m\to\iy$,
	\begin{equation}\label{11.23.13}
		\begin{split}
			\vp_\de^m &\ra \wt\vp_\de  \text{ weakly-* in } L^\iy(0,T; L^2(\Om_\de)),\\
			\vp_\de^m &\ra \wt\vp_\de  \text{ weakly in } L^2(0,T; H_{\Ga_\de}^1(\Om_\de;A_\de)), \\
			\pt_t\vp_\de^m &\ra \pt_t\wt \vp_\de \text{ weakly in } L^2(0,T; H_{\Ga_\de}^{-1}(\Om_\de;A_\de)). 
		\end{split}
	\end{equation}
	Moreover, by standard Sobolev embedding in Bochner spaces (see, e.g., \cite[Section 5.9.2, Theorem 3, p.~303]{Evans}), we have 
	\begin{equation}\label{11.23.14}
		\wt\vp_\de\in C([0,T]; L^2(\Om_\de)). 
	\end{equation}

\textit{Step 4}. We prove that $\wt\vp_\de$ is a weak solution of \eqref{11.23.2} and satisfies \eqref{11.23.4}.

Fix an integer $N\in\N^*$ and choose a test function $v\in C^1([0,T]; H_{\Ga_\de}^1(\Om_\de;A_\de))$ of the form
\begin{equation*}
	v(x,t)=\sum_{n=1}^N d_n(t)\Phi_n^\de(x), \quad d_n\in C^1([0,T]), \ n=1,\cdots, N.
\end{equation*}
For any $m\geq N$, multiplying \eqref{11.23.9} by $d_n(t)$, summing over $n=1,\cdots, N$, and integrating over $(0,T)$, we obtain
\begin{equation}\label{11.23.15}
	\int_0^T \lg \pt_t\vp_\de^m, v\rg_{H_{\Ga_\de}^{-1}(\Om_\de;A_\de), H_{\Ga_\de}^1(\Om_\de;A_\de)}\,\df t + \int_0^T (\vp_\de^m, v)_{H_{\Ga_\de}^1(\Om_\de;A_\de)}\,\df t = \int_0^T (f_\de^m, v)_{L^2(\Om_\de)}\,\df t.
\end{equation}
Passing to the limit as $m\to\iy$ along the subsequence in \eqref{11.23.13}, we get
\begin{equation}\label{11.23.16}
	\int_0^T \lg \pt_t\wt\vp_\de, v\rg_{H_{\Ga_\de}^{-1}(\Om_\de;A_\de), H_{\Ga_\de}^1(\Om_\de;A_\de)}\,\df t + \int_0^T (\wt\vp_\de, v)_{H_{\Ga_\de}^1(\Om_\de;A_\de)}\,\df t = \int_0^T (f_\de, v)_{L^2(\Om_\de)}\,\df t.
\end{equation}
Since functions of the form $\sum_{n=1}^N d_n(t)\Phi_n^\de(x)$ are dense in $L^2(0,T; H_{\Ga_\de}^1(\Om_\de;A_\de))$ (as well as in $C^1([0,T]; H_{\Ga_\de}^1(\Om_\de;A_\de))$), equation \eqref{11.23.16} holds for all $v\in L^2(0,T; H_{\Ga_\de}^1(\Om_\de;A_\de))$.

To verify the initial condition $\wt\vp_\de(0)=\vp_\de^0$, take any $v\in C^1([0,T]; H_{\Ga_\de}^1(\Om_\de;A_\de))$ with $v(T)=0$. Integrating by parts in time in \eqref{11.23.15} and \eqref{11.23.16} yields
\begin{equation*}
	-\int_0^T (\vp_\de^m, \pt_t v)_{L^2(\Om_\de)}\,\df t + \int_0^T (\vp_\de^m, v)_{H_{\Ga_\de}^1(\Om_\de;A_\de)}\,\df t = \int_0^T (f_\de^m, v)_{L^2(\Om_\de)}\,\df t + (\vp_\de^m(0), v(0))_{L^2(\Om_\de)}
\end{equation*}
and
\begin{equation*}
	-\int_0^T (\wt\vp_\de, \pt_t v)_{L^2(\Om_\de)}\,\df t + \int_0^T (\wt\vp_\de, v)_{H_{\Ga_\de}^1(\Om_\de;A_\de)}\,\df t = \int_0^T (f_\de, v)_{L^2(\Om_\de)}\,\df t + (\wt\vp_\de(0), v(0))_{L^2(\Om_\de)}.
\end{equation*}
Comparing these two identities as $m\to\iy$ and noting that $\vp_\de^m(0)\to \vp_\de^0$ strongly in $L^2(\Om_\de)$, we find
\begin{equation*}
	(\wt\vp_\de(0)-\vp_\de^0, v(0))_{L^2(\Om_\de)}=0
\end{equation*}
for all $v(0)\in H_{\Ga_\de}^1(\Om_\de;A_\de)$. By the density of $H_{\Ga_\de}^1(\Om_\de;A_\de)$ in $L^2(\Om_\de)$, we conclude that $\wt\vp_\de(0)=\vp_\de^0$ in $L^2(\Om_\de)$. Thus, $\wt\vp_\de$ is a weak solution to \eqref{11.23.2}. Furthermore, by the weak lower semicontinuity of norms together with \eqref{11.23.11}-\eqref{11.23.13}, $\wt\vp_\de$ satisfies \eqref{11.23.4}.

\textit{Step 5}. Uniqueness of the weak solution.

Suppose that $\vp_{\de,1}$ and $\vp_{\de,2}$ are two weak solutions to \eqref{11.23.2} associated with the same data $(\vp_\de^0,f_\de)$, and set $\psi_\de:=\vp_{\de,1}-\vp_{\de,2}$. Then $\psi_\de$ is a weak solution to \eqref{11.23.2} with zero initial data $\psi_\de(0)=0$ and zero source term $f_\de=0$. Taking $v=\psi_\de$ as the test function (or testing against $\psi_\de$ on $(0,t)$ for any $t\in (0,T]$), we obtain
\begin{equation*}
	\f{1}{2}\|\psi_\de(t)\|_{L^2(\Om_\de)}^2 + \int_0^t \|\psi_\de(s)\|_{H_{\Ga_\de}^1(\Om_\de;A_\de)}^2\,\df s = 0,
\end{equation*}
which implies $\psi_\de\equiv 0$ in $Q_\de$. Hence, the weak solution $\vp_\de:=\wt\vp_\de$ is unique (and consequently, the whole sequence $\{\vp_\de^m\}_{m=1}^\iy$ converges to $\vp_\de$ in the sense of \eqref{11.23.13}).

\textit{Step 6}. Higher-order regularity and estimate \eqref{11.23.5}.

Now assume further that $\vp_\de^0\in H_{\Ga_\de}^1(\Om_\de;A_\de)$ and $f_\de\in H^1(0,T; L^2(\Om_\de))$. Multiplying both sides of \eqref{11.23.9} by $\f{\df}{\df t}\vp_{\de,n}^m(t)$ and summing over $n=1,\cdots, m$, we get
\begin{equation*}
	\|\pt_t\vp_\de^m\|_{L^2(\Om_\de)}^2 + \f{1}{2}\f{\df}{\df t}\|\vp_\de^m\|_{H_{\Ga_\de}^1(\Om_\de;A_\de)}^2 = (f_\de^m, \pt_t\vp_\de^m)_{L^2(\Om_\de)} \leq \f{1}{2}\|f_\de^m\|_{L^2(\Om_\de)}^2 + \f{1}{2}\|\pt_t\vp_\de^m\|_{L^2(\Om_\de)}^2,
\end{equation*}
which yields
\begin{equation}\label{11.23.17}
	\|\pt_t\vp_\de^m\|_{L^2(\Om_\de)}^2 + \f{\df}{\df t}\|\vp_\de^m\|_{H_{\Ga_\de}^1(\Om_\de;A_\de)}^2 \leq \|f_\de^m\|_{L^2(\Om_\de)}^2.
\end{equation}
Integrating \eqref{11.23.17} over $(0,t)$ for any $t\in (0,T]$ and using the orthogonal projection properties $\|\vp_\de^m(0)\|_{H_{\Ga_\de}^1(\Om_\de;A_\de)}\leq \|\vp_\de^0\|_{H_{\Ga_\de}^1(\Om_\de;A_\de)}$ and $\|f_\de^m\|_{L^2(Q_\de)}\leq \|f_\de\|_{L^2(Q_\de)}\leq \|f_\de\|_{H^1(0,T; L^2(\Om_\de))}$, we obtain
\begin{equation*}
	\|\vp_\de^m\|_{L^\iy(0,T; H_{\Ga_\de}^1(\Om_\de;A_\de))}^2 + \|\pt_t\vp_\de^m\|_{L^2(Q_\de)}^2 \leq 2\left(\|\vp_\de^0\|_{H_{\Ga_\de}^1(\Om_\de;A_\de)}^2 + \|f_\de\|_{H^1(0,T; L^2(\Om_\de))}^2\right).
\end{equation*}
Passing to the limit as $m\to\iy$ by weak lower semicontinuity, we deduce that $\vp_\de\in L^\iy(0,T; H_{\Ga_\de}^1(\Om_\de;A_\de))$, $\pt_t\vp_\de\in L^2(Q_\de)$, and
\begin{equation}\label{11.23.18}
	\|\vp_\de\|_{L^\iy(0,T; H_{\Ga_\de}^1(\Om_\de;A_\de))} + \|\pt_t\vp_\de\|_{L^2(Q_\de)} \leq C\left(\|\vp_\de^0\|_{H_{\Ga_\de}^1(\Om_\de;A_\de)} + \|f_\de\|_{H^1(0,T; L^2(\Om_\de))}\right).
\end{equation}

Finally, from the equation $\mcA_\de \vp_\de = f_\de - \pt_t\vp_\de$ in $Q_\de$, it follows that $\mcA_\de\vp_\de\in L^2(Q_\de)$ and
\begin{equation*}
	\|\mcA_\de\vp_\de\|_{L^2(Q_\de)} \leq \|f_\de\|_{L^2(Q_\de)} + \|\pt_t\vp_\de\|_{L^2(Q_\de)} \leq C\left(\|\vp_\de^0\|_{H_{\Ga_\de}^1(\Om_\de;A_\de)} + \|f_\de\|_{H^1(0,T; L^2(\Om_\de))}\right).
\end{equation*}
Fix the smooth domain $D:=\Omega_{\delta_0}$. By nesting,
$D\subset\Omega_\delta$ for every $\delta\in[0,\delta_0]$.
Choose $\zeta=\zeta(x_2)\in C^\infty(\mathbb R)$ with
$0\leq\zeta\leq1$, $\zeta=0$ for $x_2\leq5\delta_0/2$,
and $\zeta=1$ for $x_2\geq7\delta_0/2$.
Put $\psi_\delta=\zeta\varphi_\delta|_D$.
The artificial boundary of $D$ is contained in
$\{x_2\leq2\delta_0\}$, where $\zeta$ vanishes in a
neighborhood. On the remaining boundary,
$\varphi_\delta$ has zero Dirichlet trace. Hence
$\psi_\delta\in L^2(0,T;H_0^1(D))$.

For almost every $t\in(0,T)$, the product rule gives
\[
\begin{cases}
	-\operatorname{div}(A\nabla\psi_\delta)=G_\delta
	&\text{in }D,\\
	\psi_\delta=0&\text{on }\partial D,
\end{cases}
\]
where
\[
\begin{split}
	G_\delta
	&=\zeta(f_\delta-\partial_t\varphi_\delta)
	-2\nabla\zeta\cdot A\nabla\varphi_\delta
	+\varphi_\delta\mathcal A\zeta\\
	&=\zeta\mathcal A_\delta\varphi_\delta
	-2x_2^\alpha\zeta'\partial_{x_2}\varphi_\delta
	-(x_2^\alpha\zeta')'\varphi_\delta.
\end{split}
\]
Here the restrictions to $D$ are understood. On the fixed
$C^3$ domain $D$, the matrix $A$ is uniformly elliptic and
has bounded $C^1$ coefficients. Therefore the Dirichlet
boundary $H^2$ estimate and the Poincar\'e inequality imply
\[
\begin{split}
	\|\psi_\delta(t)\|_{H^2(D)}
	&\leq C\big(\|G_\delta(t)\|_{L^2(D)}
	+\|\psi_\delta(t)\|_{L^2(D)}\big)\\
	&\leq C\big(
	\|\mathcal A_\delta\varphi_\delta(t)\|_{L^2(\Omega_\delta)}
	+\|\varphi_\delta(t)\|_{H_{\Gamma_\delta}^1
		(\Omega_\delta;A_\delta)}\big),
\end{split}
\]
with $C$ independent of $\delta$.
Since $\Omega_{4\delta_0}\subset D\cap\{x_2>4\delta_0\}$,
$\psi_\delta=\varphi_\delta$ on $\Omega_{4\delta_0}$.
Integrating in time and using the graph norm on
$D(\mathcal A_\delta)$, we obtain
\[
\begin{split}
	&\|\varphi_\delta\|_{L^2(0,T;D(\mathcal A_\delta))}
	+\|\varphi_\delta\|_{L^2(0,T;H^2(\Omega_{4\delta_0}))}\\
	&\quad\leq C\big(
	\|\mathcal A_\delta\varphi_\delta\|_{L^2(Q_\delta)}
	+\|\varphi_\delta\|_{L^2(0,T;H_{\Gamma_\delta}^1
		(\Omega_\delta;A_\delta))}\big).
\end{split}
\]
Combining this with \eqref{11.23.4} and \eqref{11.23.18}
proves \eqref{11.23.5}. The constant may depend on
$\alpha$, $M$, the fixed domain $D=\Omega_{\delta_0}$,
and $T$, but is independent of $\delta\in[0,\delta_0]$. This completes the proof of Theorem \ref{11.23.T2}. 
\end{proof}

\begin{proof}[Proof of Proposition \ref{11.23.P1}]
	First, we prove   the necessity. 
	
If $\vp_\de$ is a weak solution to \eqref{11.23.2} associated with $(\vp_\de^0, f_\de)$, then for any $\psi\in C^\iy(\ol Q_\de)$ with $\dist(\supp \psi(\cdot,t),\ol{\Ga_\de})>0$ for all $t\in [0,T]$ and $\psi(\cdot,T)=0$, we have $\psi\in C^1([0,T]; H_{\Ga_\de}^1(\Om_\de;A_\de))$. Testing \eqref{11.23.2} with $\psi$, integrating over $(0,T)$, and applying the integration-by-parts formula in time together with $\vp_\de(\cdot,0)=\vp_\de^0$ and $\psi(\cdot,T)=0$, we immediately obtain \eqref{11.23.26}.
	
	Second, we prove the sufficiency. 
	
	Conversely, suppose that $\vp_\de\in C([0,T]; L^2(\Om_\de))\cap L^2(0,T; H_{\Ga_\de}^1(\Om_\de;A_\de))$ satisfies \eqref{11.23.26}. For any $\eta\in C_0^\iy(0,T)$ and $w\in C_{\Ga_\de}^\iy(\ol{\Om_\de})$, taking $\psi(x,t)=\eta(t)w(x)$ in \eqref{11.23.26} yields
	\begin{equation}\label{11.23.27}
		\begin{split} 
		&-\int_0^T (\vp_\de(\cdot,t), w)_{L^2(\Om_\de)} \eta'(t)\,\df t \\
		&= \int_0^T \left[(f_\de(\cdot,t), w)_{L^2(\Om_\de)} - \int_{\Om_\de} \nabla\vp_\de(\cdot,t)\cdot A_\de\nabla w\,\df x\right] \eta(t)\,\df t.
		\end{split} 
	\end{equation}
	Since $C_{\Ga_\de}^\iy(\ol{\Om_\de})$ is dense in $H_{\Ga_\de}^1(\Om_\de;A_\de)$, identity \eqref{11.23.27} holds for all $w\in H_{\Ga_\de}^1(\Om_\de;A_\de)$. By the Cauchy-Schwarz and Poincar\'e inequalities, the right-hand side defines a functional in $L^2(0,T; H_{\Ga_\de}^{-1}(\Om_\de;A_\de))$. Hence, by the definition of distributional time derivatives, $\pt_t\vp_\de\in L^2(0,T; H_{\Ga_\de}^{-1}(\Om_\de;A_\de))$ and
	\begin{equation}\label{11.23.28}
		\lg \pt_t\vp_\de(t), w\rg_{H_{\Ga_\de}^{-1}(\Om_\de;A_\de), H_{\Ga_\de}^1(\Om_\de;A_\de)} + \int_{\Om_\de} \nabla\vp_\de(\cdot,t)\cdot A_\de\nabla w\,\df x = (f_\de(\cdot,t), w)_{L^2(\Om_\de)}
	\end{equation}
	for a.e. $t\in (0,T)$ and all $w\in H_{\Ga_\de}^1(\Om_\de;A_\de)$.
	
	Finally, for any $w\in C_{\Ga_\de}^\iy(\ol{\Om_\de})$ and $\zeta\in C^\iy([0,T])$ with $\zeta(0)=1$ and $\zeta(T)=0$, testing \eqref{11.23.28} with $\psi(x,t)=\zeta(t)w(x)$, integrating over $(0,T)$, integrating by parts in time, and comparing with \eqref{11.23.26} gives
	\begin{equation*}
		(\vp_\de(\cdot,0) - \vp_\de^0, w)_{L^2(\Om_\de)} = 0 \quad \text{for all } w\in C_{\Ga_\de}^\iy(\ol{\Om_\de}).
	\end{equation*}
	By the density of $C_{\Ga_\de}^\iy(\ol{\Om_\de})$ in $L^2(\Om_\de)$, we obtain $\vp_\de(\cdot,0)=\vp_\de^0$ in $L^2(\Om_\de)$. Thus, $\vp_\de$ is a weak solution to \eqref{11.23.2} associated with $(\vp_\de^0, f_\de)$. We prove the proposition. 
\end{proof}

\section*{Conflict of interest}
The authors declare that they have no conflict of interest.

\section*{Data availability}
No external datasets were used in this article. All results are established
analytically.

\end{document}